%% file: varsphere.tex
\documentclass[a4paper,11pt,oneside,fleqn]{article}

\usepackage[mathscr]{eucal}

\usepackage{amsmath,amssymb}

\usepackage{diagrams}

\usepackage{graphicx}
\usepackage{caption}
\usepackage{subcaption}
\usepackage{color}
\usepackage{standalone}
\usepackage{framed}
\usepackage{tikz,pgf}
\usetikzlibrary{patterns}
\usetikzlibrary{calc}
\usepackage[super]{nth}
\usepackage{cancel}
\usepackage{soul}
\usepackage{pdfpages}
\usepackage{hyperref}
\hypersetup{colorlinks, linkcolor=blue, citecolor=blue, urlcolor=blue}

\definecolor{Cyan}{rgb}{0.0, 0.55, 0.55}

\usepackage{lineno}

\newcommand{\changed}[1]{\textcolor{black}{#1}}

\newcommand{\const}{\mathop{\rm constant}\nolimits}

\DeclareMathOperator\Div{Div}
\DeclareMathOperator\Curl{Curl}

\newtheorem{theorem}{Theorem}

\newtheorem{proposition}[theorem]{Proposition}

\title{Variational integrator for the rotating shallow-water equations on the sphere}

\begin{document}

	\par\noindent {\LARGE\bf
	Variational integrator for the rotating shallow-water equations on the sphere
	\par}

{\vspace{4mm}\par\noindent {\bf R\"udiger Brecht$^\dag$, Werner Bauer$^\ddag$, Alexander Bihlo$^\dag$, Fran\c cois Gay-Balmaz$^\S$ and Scott MacLachlan$^\dag$ 
	} \par\vspace{2mm}\par}

{\vspace{2mm}\par\noindent {\it
		$^{\dag}$~Department of Mathematics and Statistics, Memorial University of Newfoundland,\\ St.\ John's (NL) A1C 5S7, Canada
}}
{\vspace{2mm}\par\noindent {\it
		$^{\ddag}$~Imperial College London, Department of Mathematics, 180 Queen’s Gate, London SW7 2AZ, United Kingdom.
}}
{\vspace{2mm}\par\noindent {\it
		$^{\S}$~\'Ecole Normale Sup\'erieure/CNRS, Laboratoire de M\'et\'eorologie Dynamique, Paris, France.
}}

{\vspace{2mm}\par\noindent {\it
		\textup{E-mail:} rbrecht@mun.ca, w.bauer@imperial.ac.uk, abihlo@mun.ca, gaybalma@lmd.ens.fr, smaclachlan@mun.ca 
	}\par}

\vspace{4mm}\par\noindent\hspace*{10mm}\parbox{140mm}{\small
	We develop a variational integrator for the shallow-water equations on a rotating sphere. The variational integrator is built around a discretization of the continuous Euler--Poincar\'{e} reduction framework for Eulerian hydrodynamics. We describe the discretization of the continuous Euler--Poincar\'{e} equations on arbitrary simplicial meshes. Standard numerical tests are carried out to verify the accuracy and the excellent conservational properties of the discrete variational integrator. 
}\par\vspace{4mm}

	\section{Introduction}
	
	Geometric numerical integration is the branch of numerical analysis devoted to the development of discretization schemes that preserve important geometric properties of differential equations. Examples of geometric properties of practical interest include symplectic forms and Hamiltonian formulations, Lie group symmetries and conservation laws, volume forms, variational formulations, maximum principles and blow-up properties. Numerical integrators that discretely preserve one or more of the aforementioned geometric properties are presented e.g.\ in~\cite{baue17a,baue17b,beck15a,bihl17b,blan16a,budd01a,hair06Ay,leim04Ay,mars01a,pavl11a,sanz94a,wan16a}.
	
	A main motivation behind the development of integrators capable of preserving geometric properties of differential equations is their, in general, superior long-term behaviour. Preserving geometric properties can guarantee arbitrarily long-term stability, consistency in statistical properties and the prevention of systematic drift in stationary or periodic solutions, see e.g.~\cite{hair06Ay,leim04Ay,wan16b} for further details. 
	
	Recent years have seen an increased interest in geometric numerical integration for models of atmospheric dynamics. This is natural since long time integrations and the accurate representation of the statistical properties of these models lie at the heart of climate prediction and turbulence modeling. A particular model that has received considerable attention is the rotating shallow-water (RSW) equations, both in the plane and on the sphere.
	
	Energy and enstrophy preserving integrators for the shallow-water equations were developed as early as in the seminal paper of \cite{arak81a}. There it was recognized that preserving energy and enstrophy in a finite difference discretization of the shallow-water equations is crucial to guaranteeing the numerical stability of typical flow regimes. In recent years, considerable effort was devoted to the development of structure-preserving integrators on general structured and unstructured grids, see e.g.~\cite{baue17a,bihl12By,eldr17a,fran03a,ring02a,thub12a,thub15a}. 
	
	In~\cite{baue17a}, a variational integrator for the rotating shallow-water equations in the plane was proposed. Variational integrators rest on first discretizing the continuous variational principle underlying the governing equations of interest, and then deriving the numerical scheme as a discrete system of Euler--Lagrange equations~\cite{mars01a}. Variational integrators possess a number of desirable properties, including compatibility with a discrete form of Noether's theorem that guarantees the exact numerical preservation of those conserved quantities related to the variational symmetries of the discretized governing equations, as well as stability for exponentially long time periods~\cite{hair06Ay,mars01a}. 
	
	While most of the work on variational integration was devoted to ordinary differential equations, recent years have seen an increased interest in the partial differential equation case, see e.g.~\cite{pavl11a} and~\cite{baue17a,baue17b,DeGaGBZe2014} for some applications of the variational methodology to important models of geophysical fluid dynamics. Variational integrators for the partial differential equations of fluid dynamics are designed by replacing the continuous configuration space of the model equation, represented as an appropriate infinite-dimensional Lie group, by a suitable finite dimensional matrix Lie group on which the variational principle can be applied in both its Lagrangian and Eulerian versions, thanks to an application of the Euler--Poincar\'{e} reduction theorem. The purpose of the present paper is to extend the variational integrator proposed in~\cite{baue17a} for the shallow-water equations in the plane to the shallow-water equations on the rotating sphere.
	
	The further organization of the paper is as follows. In Section~\ref{sec:VariationalPrinciple}, we give a brief summary of the variational description of the shallow-water equations using the Euler--Poincar\'{e} formulation.
	Section~\ref{sec:DiscreteVariationalPrinciple} is devoted to the description of the discretization of the continuous Euler--Poincar\'{e} formulation, originally presented in~\cite{baue17a}, and the representation of the variational integrator on the icosahedral mesh geometry used to approximate the sphere.
	In Section~\ref{sec_notations}, we verify the consistency of the corresponding approximations of the standard differential operators. Test cases and numerical benchmarks showcasing the behaviour of the variational integrator for the shallow-water equations are given in Sections~\ref{sec:NumericalSimulationsVariationalIntegrator} and \ref{sec:LongTermSimulations}. The conclusions and thoughts for future research within this field of geometric numerical integration are found in Section~\ref{sec:ConclusionsVariationalIntegrator}.

	\section{Variational principle for the rotating shallow-water equations}\label{sec:VariationalPrinciple}

	In absence of irreversible processes, the equations of motion of fluid dynamics can be derived via the Hamilton principle
	\[
	\delta \int_0^T L (\varphi, \dot\varphi)dt=0,
	\]
	\changed{with respect to variations $\delta \varphi$ vanishing at $t=0$ and $t=T$. Here $L$ is the Lagrangian function of the fluid, given by the kinetic energy minus the internal energy, and expressed in terms of the Lagrangian fluid trajectory $\varphi$ and Lagrangian fluid velocity $\dot\varphi$.} Following this point of view, the configuration space for fluid dynamics, away from shocks, is the group $\operatorname{Diff}(D)$ of diffeomorphisms $\varphi$ of the fluid domain $D$.
	While in the Lagrangian (or material) description this principle is a straightforward extension of the Hamilton
	principle of classical mechanics, in the Eulerian (or spatial) description the variational principle is more involved, since it uses constrained variations. It can be rigorously justified by applying the process of Euler--Poincar\'e reduction \cite{holm98Ay}, which directly gives the general form of the Eulerian variational principle induced by the Hamilton principle of fluid mechanics. We refer to \cite{holm99Ay} for an application to the equations of geophysical fluid dynamics (GFD). Since this principle is central to the derivation of the numerical scheme, we quickly review it below for the shallow-water case. \changed{A more detailled review of Euler--Poincar\'e reduction and its application to the shallow water equations as well as to simpler examples is presented in \S\ref{Appendix_EP_principle} in the Appendix.}

	Let us consider the rotating shallow-water (RSW) dynamics on the sphere $\mathcal{S}$ of radius $R$. The sphere is naturally endowed with a Riemannian metric $\gamma$ induced from the standard Euclidian metric on $\mathbb{R}^3$ and with a volume $d\sigma$ associated to $\gamma$, \changed{see \S\ref{Appendix_forms}}. In terms of latitude ($\theta$) and longitude ($\lambda$), we have $\gamma= R^2 d\theta^2 + R^2\cos \theta d \lambda^2$ and $d\sigma = R^2 \cos \theta d\theta\wedge d\lambda$, but our approach is geometrically intrinsic, i.e., independent of any choice of coordinates on $\mathcal{S}$. 
	In the Eulerian description, the variables are the fluid velocity $\mathbf{u}$ and the fluid depth $h$, defined in terms of the Lagrangian variables $\varphi$ and~$\dot \varphi$, with $\varphi \in\operatorname{Diff}(\mathcal{S})$, as	\begin{equation}\label{def_u_h}
	\mathbf{u}=\dot\varphi\circ\varphi^{-1}\quad\text{and}\quad h= h_0\bullet \varphi^{-1}:=(h_0\circ \varphi^{-1}) J\varphi^{-1},
	\end{equation}
	where $h_0$ is the initial fluid depth and $J\varphi^{-1}$ is the Jacobian of the diffeomorphism $\varphi^{-1}$ with respect to the metric. The second relation in \eqref{def_u_h} is the \emph{natural action} of diffeomorphisms on densities, that we have denoted using $\bullet$, \changed{see also \S\ref{Appendix_forms} for more details}.
	
	The Lagrangian for rotating shallow-water fluids in Eulerian coordinates is given by
	\begin{equation}\label{Lagrangian_barotropic}
	\ell( \mathbf{u} , h)= \int_{\mathcal{S}} \Big[ \frac{1}{2} h \, \gamma(\mathbf{u} ,\mathbf{u}) + h \,\gamma(\mathbf{R},\mathbf{u})
	- \frac{1}{2} g(h + B) ^2 \Big]d \sigma,
	\end{equation}
	where $B$ is the bottom topography, $h+B$ describes the free surface elevation of the fluid, $g$ is the gravitational acceleration and $ \mathbf{R} $ is the vector potential of the angular velocity of the Earth. \changed{We recall that $\gamma$ is the Riemannian metric on $\mathcal{S}$ induced from the Euclidean metric on $\mathbb{R}^3$}.
	
	Given this Lagrangian, the equations of motion follow from the Euler--Poincar\'e variational principle given by
	\begin{equation}\label{EP_principle}
	\delta\int_0^T\ell(\mathbf{u},h)dt=0,
	\end{equation}
	with respect to constrained variations of the form $\delta\mathbf{u}=\partial_t\mathbf{v}+ [\mathbf{u},\mathbf{v}]$ and $\delta h=-\operatorname{div}(h\mathbf{v})$, where $\mathbf{v}$ is an arbitrary time dependent vector field on $\mathcal{S}$, vanishing at the endpoints $t=0$ and $t=T$, and where $\operatorname{div}$ denotes the divergence on the sphere, associated to the metric, \changed{see \S\ref{Appendix_forms}}. 
	The form of the constrained variations is obtained by using the relations \eqref{def_u_h} and computing the variations $\delta\mathbf{u}$ and $\delta h$ induced by free variations $\delta\varphi$ vanishing at $t=0$ and $t=T$. This principle yields the Euler--Poincar\'e equations in the general form
	\begin{equation}\label{EP_general} 
	\partial _t \frac{\delta \ell}{\delta \mathbf{u} }+  \mathsf{L} _\mathbf{u} \frac{\delta \ell}{\delta \mathbf{u} } =h\,\mathbf{d}  \frac{\delta \ell}{\delta h},
	\end{equation}
	where the second term denotes the Lie derivative of the fluid momentum density 
	(a one-form density) along the vector field $\mathbf{u}$, explicitly given by $\mathsf{L}_\mathbf{u} \alpha=\mathbf{i}_\mathbf{u} \mathbf{d}\alpha+ \mathbf{d}\mathbf{i}_\mathbf{u}\alpha+\alpha\operatorname{div}\mathbf{u}$, with $\mathbf{d}$ being the exterior derivative and \changed{$\mathbf{i}_\mathbf{u}\omega$ the inner product of the $k$-form $\omega$ with the vector field $\mathbf{u}$. We refer to \S\ref{Appendix_forms} for a review of exterior derivatives and Lie derivatives and to \S\ref{Appendix_EP_principle} for a detailed derivation of \eqref{EP_general}}. The functional derivative $\frac{\delta \ell}{\delta \mathbf{u} }$, resp. $\frac{\delta \ell}{\delta h }$, is the one-form, resp., the function, defined by
	\[
	\int_\mathcal{S} \frac{\delta \ell}{\delta \mathbf{u} }\cdot   \mathbf{v}\,d\sigma= \left.\frac{d}{d\varepsilon}\right|_{\varepsilon=0}\ell(\mathbf{u}+\varepsilon\mathbf{v},h),\quad \text{resp.}\quad  \int_\mathcal{S} \frac{\delta \ell}{\delta h } v \,d\sigma= \left.\frac{d}{d\varepsilon}\right|_{\varepsilon=0}\ell(\mathbf{u},h+\varepsilon v)
	\]
	for any vector field $\bf v$ or function $v$.
	Equation \eqref{EP_general} is accompanied with the advection equation $\partial _t h + \operatorname{div}(h \mathbf{u} )=0$ for the fluid depth which, in the Euler--Poincar\'e approach, follows from the second relation in \eqref{def_u_h}.
	%
	
	For the RSW Lagrangian \eqref{Lagrangian_barotropic}, we have
	\begin{equation}\label{FD_RSW}
	\frac{\delta \ell}{\delta \mathbf{u} }=h ( \mathbf{u} ^\flat+ \mathbf{R}^\flat)\quad\text{and}\quad \frac{\delta \ell}{\delta h}= \frac{1}{2} \gamma(\mathbf{u}, \mathbf{u})+ \gamma(\mathbf{R}, \mathbf{u})- g(h+B),
	\end{equation}
	\changed{where $\mathbf{u}^\flat$ and $\mathbf{R}^\flat$ are the one-forms associated to the vector fields $\mathbf{u}$ and $\mathbf{R}$ via the flat operator $\flat$ of the Riemannian metric $\gamma$, see \S\ref{Appendix_forms}}. Using these expressions, the Euler--Poincar\'e equations \eqref{EP_general}
	lead to the momentum RSW equations, written in the space of one-forms as
	\begin{equation}\label{one_form_RSW} 
	h \,\partial _t \mathbf{u} ^\flat + \mathbf{i} _{ h \mathbf{u} } \mathbf{d} ( \mathbf{u}^\flat 
	+ \mathbf{R}^\flat )=- h\,\mathbf{d} \Big(  \frac{1}{2} \gamma( \mathbf{u} , \mathbf{u}) + g(h+B) \Big),
	\end{equation}
	see \cite{baue17a} \changed{and \S\ref{Appendix_EP_principle}} for details.
	It is this expression of the RSW equations that appears in a discretized form in the variational discretization later in \eqref{simplicial_scheme}.

	\section{Discrete variational principle for the RSW equations}
	\label{sec:DiscreteVariationalPrinciple}

	The variational discretization of the rotating shallow-water equations 
	mimics the continuous variational method; in particular, each step of the continuous theory is translated
	to the discrete level. We provide here a review of the discrete Euler--Poincar\'e theory for the RSW, and refer the 
	reader to \cite{baue17a} for full details.

	\subsection{Definition of the appropriate discrete configuration space for RSW}
	\changed{\paragraph{Discrete diffeomorphism group.}
		The discretization procedure starts with the choice of a discrete version of the configuration group $\operatorname{Diff}(\mathcal{S})$ of a shallow-water fluid. Following  \cite{pavl11a}, \cite{baue17a}, this group is obtained by first discretizing the space of functions $\mathcal{F}(\mathcal{S})$ on which $\operatorname{Diff}(\mathcal{S})$ acts by composition on the right, and then identifying a finite dimensional group acting by matrix multiplication on the finite dimensional space of discrete functions, while preserving some properties of the action by diffeomorphisms. In the compressible case, we shall only retain the property that constant functions are preserved under composition by a diffeomorphism. Given a mesh $\mathbb{M}$ of the sphere and choosing as discrete functions the space $\mathbb{R}^N$ of piecewise constant functions on~$\mathbb{M}$, this results in the discrete diffeomorphism group
		\begin{equation}\label{DD_all}
		\mathsf{D}(\mathbb{M})=\left\{q\in \operatorname{GL}(N)^+\mid q\cdot \boldsymbol{1}=\boldsymbol{1}\right\}
		\end{equation}
		of dimension $ N ^2 -N$, where $N$ is the number of cells in $\mathbb{M}$, the vector $\boldsymbol{1}\in\mathbb{R}^N$ is defined by $\boldsymbol{1}=(1,...,1)^\mathsf{T}$, and $\operatorname{GL}(N)^+=\{q\in \operatorname{Mat}(N)\mid \operatorname{det}q>0\}$ is the general linear group of orientation-preserving $N\times N$ matrices.
		The condition $q\cdot  \boldsymbol{1}= \boldsymbol{1}$ encodes, at the discrete level, the fact that constant functions are preserved under composition by a diffeomorphism. The action of the group $\mathsf{D}(\mathbb{M})$ by matrix multiplication on discrete functions $F\in \mathbb{R}^N$ is denoted as				\[
		F \in \mathbb{R}^N\mapsto qF= F\circ q^{-1} \in\mathbb{R}^N,\quad q\in \mathsf{D}(\mathbb{M}),
		\]
		where the suggestive notation  $F\circ q^{-1}$ for the multiplication of the vector $F$ by the matrix $q$ is introduced to indicate that this action is understood as a discrete version of the action of $\operatorname{Diff}(\mathcal{S})$ by composition on the space $\mathcal{F}(\mathcal{S})$ of functions on $\mathcal{S}$, namely $f\in\mathcal{F}(\mathcal{S})\mapsto f\circ\varphi^{-1}\in\mathcal{F}(\mathcal{S})$, see \cite{pavl11a}, \cite{baue17a} for details. The situation is formally illustrated by the diagram
		\[
		\begin{diagram}
		f\in \mathcal{F}(\mathcal{S}) &   &\rTo(2,2)^{\operatorname{Diff}( \mathcal{S} )\;\;\;\;\;\;}&   f\circ\varphi^{-1}\in \mathcal{F} ( \mathcal{S} )\\
		\dTo(1,1)
		&& & \dTo(1,1)
		\\
		F\in \mathbb{R} ^N & &  \rTo(2,2)^{\mathsf{D}( \mathbb{M})\;\;\;\;\; }&q F=F\circ q^{-1}\in \mathbb{R}  ^N.
		\end{diagram}
		\]			
		The Lie algebra of the Lie group $\mathsf{D}(\mathbb{M})$ is the space of row-null $N \times N$ matrices
		\begin{equation}\label{Lie_algebra} 
		\mathfrak{d}(\mathbb{M})=\{A\in \operatorname{Mat}(N)\mid A\cdot  \boldsymbol{1}=0\},
		\end{equation} 
		endowed with the Lie bracket $[A,B]=AB-BA$. This Lie algebra is a discrete version of the Lie algebra of $\operatorname{Diff}( \mathcal{S} )$ given by vector fields on $\mathcal{S}$.
		It is of particular interest for our derivations as it will allow us to formulate
		the discrete spatial Lagrangian required to derive the Euler--Poincar\'e equations 
		from variational principles.  By taking the derivative of continuous and discrete actions at the identity, we get $\left.\frac{d}{dt}\right|_{t=0}f\circ\varphi^{-1}_t=-\mathbf{d}f\cdot \mathbf{u}$ and $\left.\frac{d}{dt}\right|_{t=0}F\circ q^{-1}_t=AF$, where $\left.\frac{d}{dt}\right|_{t=0}\varphi_t=\mathbf{u}$ and $\left.\frac{d}{dt}\right|_{t=0}q_t=A$. This suggests that $AF$, with $A$ an element of the Lie algebra $\mathfrak{d}(\mathbb{M})$, may play the role of a discrete version of the derivative of $f$ in the direction of a continuous vector field $\mathbf{u}$. As we will recall below, not all $A\in \mathfrak{d}(\mathbb{M})$ can be interpreted as a discrete vector field. This induces nonholonomic constraints on the Lie algebra \eqref{Lie_algebra}, which have to be appropriately taken into account in the variational principle.
		\paragraph{Nonholonomic constraints.}} It can be shown that if a matrix $A\in \mathfrak{d}(\mathbb{M})$ approximates a vector field $\mathbf{u}$, then the matrix elements $A_{ij}$ satisfy
	\begin{equation}\label{approx_A_ij} 
	\begin{aligned} 
	A _{ij}& \simeq - \frac{1}{2 \Omega _{ii}}\int_{ D_{ij} }( \mathbf{u} \cdot \mathbf{n}_{ij}) dS, \quad \text{for all $j \in N(i)$, $j\neq i$},\\
	A _{ii} & \simeq  \frac{1}{2 \Omega _{ii}}\int_{ C _i  } (\operatorname{div} \mathbf{u} )d \mathbf{x},
	\end{aligned}
	\end{equation}
	where $N(i)$ denotes the set of all indices (including $i$) of cells sharing a face with cell $C_i$, $D_{ij}$ denotes the face common to cells $C_i$ and $C_j$ with unit normal $\mathbf{n}_{ij}$ pointing from $C_i$ to $C_j$, and~$\Omega_{ii}$ is the volume of cell $C_i$. We refer to \cite{baue17a} for the precise statement of these approximations.
	This identification imposes several constraints on the matrices in $\mathfrak{d}(\mathbb{M})$ to ensure that they represent a velocity vector field $\mathbf{u}$. First it is 
	required that fluxes are nonzero only between neighboring cells, hence we have the linear constraint
	\begin{equation}\label{constraint_1}
	\mathcal{S} = \big\{ A \in \mathfrak{d} ( \mathbb{M}  )\mid A _{ij} =0, \;\; \text{for all $j \notin N(i)$} \big \}.
	\end{equation}
	Second, we have the constraint $ \Omega _{ii} A _{ij} =- \Omega _{jj} A _{ji}$, 
	for all $j\neq i$, i.e., $A^\mathsf{T}\Omega + \Omega A$ is a diagonal matrix, with $\Omega$ being the $N\times N$ diagonal matrix with elements $\Omega_{ii}$. 
	This gives the additional linear constraint
	\begin{equation}\label{constraint_2} 
	\mathcal{R} = \big \{ A \in \mathfrak{d} ( \mathbb{M}  )\mid A^\mathsf{T} \Omega + \Omega A\;\; \text{is diagonal} \big\}.
	\end{equation}
	\changed{For the subsequent application of the variational principle, it is important to note that the constraints \eqref{constraint_1} and \eqref{constraint_2} are nonholonomic. Such constraints are taken into account by using the \textit{Euler--Poincar\'e--d'Alembert principle}, which is the nonholonomic version of the Euler--Poincar\'e principle.  As we recall in \S\ref{Appendix_EP_principle}, the application of this principle requires the choice of an appropriate dual space and duality pairing.}

	\changed{\paragraph{Duality pairing and projector.}}We recall, see \cite{pavl11a}, that in the context of the discrete diffeomorphism group, the space of discrete one-forms $\Omega_d^1(\mathbb{M})$ is identified with the space of skew-symmetric $N \times N$ matrices. The discrete version of the $L^2$-pairing between discrete one-forms and discrete vector fields is given by\begin{equation}\label{duality_pairing}
	\left\langle L, A \right\rangle = \operatorname{Tr}(L^\mathsf{T} \Omega A),\quad \text{for}\quad A,L\in \operatorname{Mat}(N).
	\end{equation}

	The application of this variational principle makes crucial use of Proposition \ref{DM} 
	of \cite{baue17a} recalled below, which identifies a projector onto the dual space to the constraint $\mathcal{R}$ with respect to the pairing \eqref{duality_pairing}. The role of this proposition will become clear below when stating the Euler--Poincar\'e--d'Alembert principle.
	
	\begin{proposition}[\cite{baue17a}]\label{DM} Given an $N\times N$ matrix $L$, we have the equivalence
		\[
		\left\langle L, A \right\rangle =0,\;\;\text{for all $A \in \mathcal{R} $}\quad \Leftrightarrow \quad \mathbf{P} (L)=0,
		\]
		where $\mathbf{P} : \operatorname{Mat}(N)\rightarrow  \Omega_d^1 ( \mathbb{M}  )$ is the projector defined by $\mathbf{P} (L):=( L- \widehat{L} )^{(A)}$, with $\widehat{L}_{ij}:= L_{ii}$ and $L^{(A)}:= \frac{1}{2} (L- L^\mathsf{T})$.
	\end{proposition}

	\subsection{Euler--Poincar\'e--d'Alembert variational principle}
	
	We shall now reproduce, at the discrete level, the variational formulation for the RSW recalled earlier in Section \ref{sec:VariationalPrinciple}, with the goal of obtaining the semi-discrete RSW equations via the Euler--Poincar\'e--d'Alembert principle.

	As a first step, we need to identify the action of the group $\mathsf{D}(\mathbb{M})$ on the variables $D\in\mathbb{R}^N$ representing the discrete fluid depth. As in the continuous case, see the second equation in \eqref{def_u_h}, this action, also denoted by $D\mapsto D\bullet q$, is dual to the action on discrete functions, namely
	\begin{equation}\label{def_advection}
	\langle D\bullet q, F\rangle_0 = \langle D, F\circ q^{-1}\rangle_0,\;\;\text{for all $F\in\mathbb{R}^N$, $q\in\mathsf{D}(\mathbb{M})$},
	\end{equation}
	with respect to the discrete $L^2$-pairing $\left\langle D, F\right\rangle_0=D^\mathsf{T}\Omega F$.  A direct application of \eqref{def_advection} shows that this action and the associated Lie algebra action are given by the formulas
	\begin{equation}\label{action_density} 
	D \bullet q=  \Omega ^{-1} q^{\mathsf{T}} \Omega D\quad\text{and}\quad  D \bullet A= \Omega ^{-1} A^{\mathsf{T}} \Omega D,
	\end{equation}
	for $q\in\mathsf{D}(\mathbb{M})$ and $A\in\mathfrak{d}(\mathbb{M})$.\color{black}

	Given a semi-discrete Lagrangian $\ell(A,D)$, the discrete version of the variational principle \eqref{EP_principle} reads as follows
	\begin{equation}\label{VP_Euler} 
	\delta \int_0^T \ell( A, D) dt=0,
	\end{equation} 
	with respect to variations $\delta A=\partial_tB+[B,A]$ and $ \delta D=- D\bullet B$, with $A , B\in \mathcal{S} \cap \mathcal{R} $ and $B(0)=B(T)=0$. As in the continuous case, the expression of the variations follows from the definition of $A$ and $D$ in terms of the Lagrangian variables $q$, $\dot q$, namely $A = \dot q q^{-1}$ and $D = D_0 \bullet q^{-1}$, which are discrete versions of \eqref{def_u_h}. This principle, which incorporates the nonholonomic constraint $\mathcal{S} \cap \mathcal{R} $ in the Euler--Poincar\'e approach, is referred to as the Euler--Poincar\'e--d'Alembert principle. \changed{Note that both $A=\dot q q ^{-1}$ and $B=\delta q q^{-1}$ have to satisfy the nonholonomic constraint $\mathcal{S} \cap \mathcal{R}$}.
	
	In the next theorem, this principle is applied to yield the general semi-discrete form of compressible fluid equations.
	
	\begin{theorem}[Discrete variational equations, \cite{baue17a}]\label{theorem_1}
		For a semi-discrete Lagrangian $\ell= \ell (A,D): \mathfrak{d} ( \mathbb{M}  ) \times \mathbb{R}^N  \rightarrow \mathbb{R}  $, the curves $A(t), D(t)$ are critical for the variational principle \eqref{VP_Euler} if and only if they satisfy the equations
		\begin{equation}\label{discrete_EP_compressible} 
		\mathbf{P} \left( \frac{d}{dt} \frac{\delta  \ell}{\delta  A}+  \Omega ^{-1} \left[ A^\mathsf{T}, \Omega  \frac{\delta  \ell}{\delta  A}\right] + D \frac{\delta \ell}{\delta D}^\mathsf{T} \right)_{ij}  =0, \quad \text{for all $i \in N(j)$},
		\end{equation}
		where $ \mathbf{P} : \operatorname{Mat}(N) \rightarrow \Omega ^1 _d ( \mathbb{M}  )$ 
		is the projection obtained in Proposition \ref{DM}. These equations are accompanied with the discrete continuity equation
		\begin{equation}\label{continuity} 
		\frac{d}{dt} D+ D\bullet A=0.
		\end{equation} 
	\end{theorem}
	
	This result follows from a direct application of \eqref{VP_Euler} by using the expression for $\delta A$ and $\delta D$ and isolating $B$ which is an arbitrary curve in $\mathcal{S}\cap\mathcal{R}$. The equations then follow by an application of Proposition \ref{DM}. 
	
	\changed{				More precisely, using the definition of the pairings $\langle\,,\rangle $ and $\langle\,,\rangle_0$, we have
		\begin{align}
		\delta \int_0^T \ell( A, D) dt&=\int_0^T \Big[ \Big\langle\frac{\delta\ell}{\delta A}, \partial_tB+[B,A]\Big\rangle+\Big\langle\frac{\delta\ell}{\delta D}, - \Omega ^{-1} B^{\mathsf{T}} \Omega D\Big\rangle_0\Big]dt\\
		&= -\int_0^T \Big \langle \frac{d}{dt} \frac{\delta  \ell}{\delta  A}+  \Omega ^{-1} \Big[ A^\mathsf{T}, \Omega  \frac{\delta  \ell}{\delta  A}\Big] + D \frac{\delta \ell}{\delta D}^\mathsf{T}, B\Big\rangle dt.
		\end{align}
		Since $B$ is arbitrary in the space $\mathcal{S} \cap \mathcal{R}$, we obtain the equations \eqref{discrete_EP_compressible} by application of Proposition \ref{DM}.
		We provide in Table \ref{table} a summary that enlightens the correspondence between the continuous and discrete objects.
	}
	
	\begin{table}[h!]
		\centering
		\begin{tabular}{|c | c |} 
			\hline
			Continuous diffeomorphisms & Discrete diffeomorphisms\\ 
			\hline
			$\operatorname{Diff}(M)\ni\varphi $ & $\mathsf{D}( \mathbb{M}  )\ni q$ \\ 
			\hline\hline
			Group action on functions &   Group action on discrete functions\\ 
			\hline
			$f \mapsto f \circ \varphi $ & $ F\circ q=F\mapsto q^{-1} F$  \\ 
			\hline\hline
				Group action on densities & Group action on discrete densities\\ 
			\hline
			$ h \mapsto h\bullet \varphi=( h  \circ \varphi)J \varphi$ & $ D\mapsto D\bullet q=\Omega^{-1} q^\mathsf{T}\Omega D$  \\ 
			\hline\hline
					Eulerian velocity and depth &  Eulerian discrete velocity and discrete depth  \\ 
			\hline
			$\mathbf{u} = \dot{ \varphi } \circ \varphi^{-1}$, \;$ h =(h _0 \circ \varphi ^{-1} ) J \varphi^{-1} $
			& $ A= \dot{ q} q^{-1}$, \;$ D =\Omega^{-1} q^{-\mathsf{T}}\Omega D_0$\\ 
			\hline\hline
				Euler-Poincar\'e principle &  Euler-Poincar\'e-d'Alembert principle\\
			\hline
			$\delta \int_0^T \ell( \mathbf{u}  ,  h ) dt=0$, &$ \delta \int_0^T \ell( A , D ) dt=0$,  \\
					$ \delta \mathbf{u} =\partial _t  \mathbf{v}+ [ \mathbf{v}, \mathbf{u} ], \;\;  \delta h = - \operatorname{div}( h\mathbf{v})$ &  $\delta A= \partial_t B+[B,A],\qquad  \delta D=- \Omega ^{-1} B^\mathsf{T} \Omega D$\\
			\qquad  & constraint: $A, B \in \mathcal{S} \cap \mathcal{R}$ \\
			\hline
		\end{tabular}
		\vspace{1em}
		\caption{Continuous and discrete objects} \label{table} 
	\end{table}
	
	For the RSW case, the discrete Lagrangian is
	\begin{equation}\label{discrete_Lagrangian} 
	\ell(A, D) = \frac{1}{2} \sum_{i,j=1}^ND_i A ^\flat_{ij} A _{ij} \Omega _{ii}+ \sum_{i,j=1}^ND_i R ^\flat_{ij} A _{ij} \Omega _{ii} - \frac{1}{2}\sum_{i=1}^N g  (D_i +  B_i)^2  \Omega _{ii},
	\end{equation} 
	see \eqref{Lagrangian_barotropic}, which requires the construction of a discrete ``flat'' operator $A\in \mathcal{S}\cap \mathcal{R}\mapsto A^\flat \in\Omega^1_d(\mathbb{M})$ associated to a given mesh, see \cite{pavl11a}. 
	
	The abstract developments made so far are valid for any kind of reasonable non-degenerate mesh
	\changed{(i.e. the mesh admits a non-degenerate circumcenter dual and 
		is member of a shape-regular and quasi-uniform mesh family, cf. \cite{baue17a})}. By choosing a fixed mesh, we will be able to express these abstract notions in concrete
	(implementable) equations.

	\subsection{Semi-discrete scheme on a 2D simplicial mesh}

	%
	
\begin{figure}[!ht]
		\centering
		\includegraphics*[width=\textwidth]{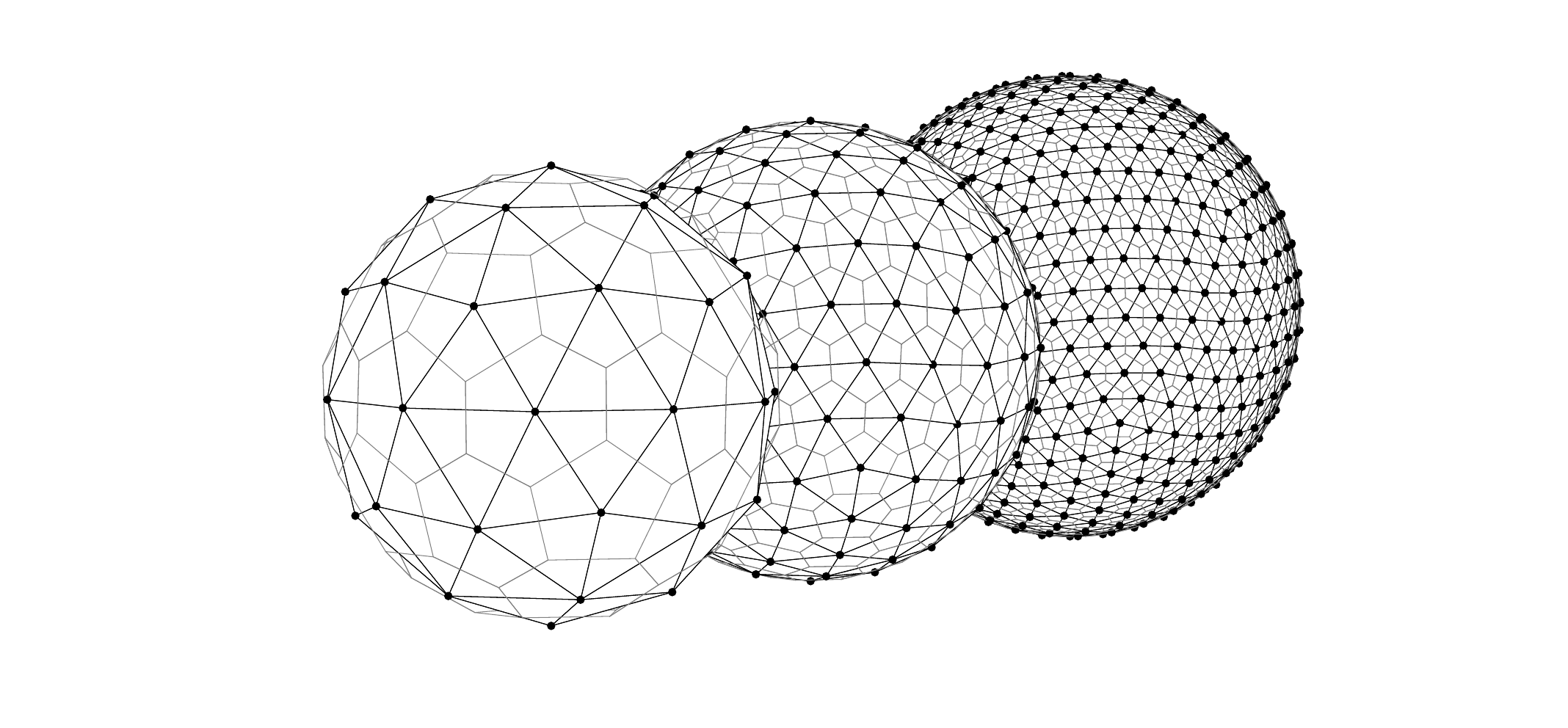}
		\caption{Icosahedral meshes with refinement levels $2$, $3$, and $4$ 
			corresponding to 80, 320, and 1280 triangular cells.}
		\label{fig:icos}
	\end{figure}

	For our implementation, we use an icosahedral grid as pictured in Fig. \ref{fig:icos}. The construction of the grid is described in \cite{Heikes_Randell_I:1995}. The icosahedron's edges are recursively bisected, and the new vertices are projected onto the unit sphere. 
	\changed{The triangles are used as the primal grid, and the circumcenter dual, consisting of pentagonal and hexagonal cells, as the dual grid. The vertex positions are optimized in the sense that the global maximum of the discrepancy between primal edge midpoints and intersection points of primal and dual edges is minimized \cite{Heikes_Randell_II:1995}. Note that, for our variational integrator, any other optimization for which primal and dual edges intersect perpendicularly could be used too; only the order of convergence of the differential operators will be affected by the optimization. 
	}
	
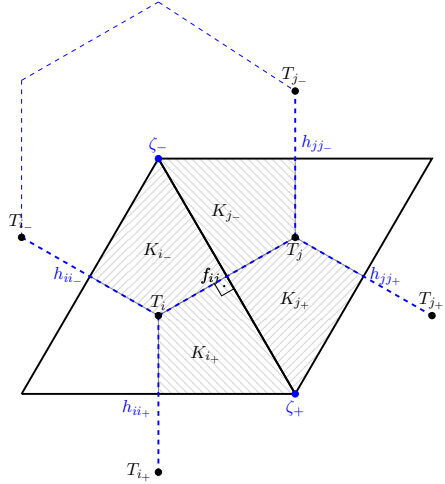
\begin{figure}[!ht]
		\centering
		\scalebox{.6}{\input{grid.tex}}
		\caption{Notation and indexing conventions for the 2D simplicial mesh.}
		\label{fig_Notation}
	\end{figure}

	In \cite{Thuburn2009}, the grid and its connectivity is described in more details.					
	Fig.~\ref{fig_Notation} shows a section of the simplicial mesh where we indicate our notation:
	
	\begin{align*} 
	f_{ij}:=  &\;\text{length of a primal edge, triangle edge located between triangle $i$ and triangle $j$};\\
	h_{ij}:=  &\;\text{length of a dual edge that connects the circumcenters of triangle $i$ and triangle $j$};\\
	\Omega_{ii}:= &\;\text{area of a primal simplex (triangle) $T_i$};\\
	K_{i}^\pm:= &\;\text{proportional area of the intersection of (triangle) $T_i$ and (hexagon/pentagon) $\zeta_\pm$},\\
	&\;\text{see Equation \eqref{equ_K}}.
	\end{align*}
	
	Using the flat operator on a 2D simplicial mesh, see \cite{pavl11a}, 
	we are able to represent the discrete one-forms $A ^\flat \in \Omega _d ^1 (\mathbb{M}  )$ in terms of the discrete vector fields $A \in\mathcal{S}\cap \mathcal{R}$ by 
	\begin{equation}\label{def_flat}
	\begin{aligned} 
	&A ^\flat _{ij}= 2 \Omega_{ii}\frac{h_{ij} }{f _{ij} }A _{ij}, \quad \text{for $j \in N(i)$},\\
	&A ^\flat _{ij}+ A ^\flat _{jk}+ A^\flat _{ki}= K_j ^{e} \big\langle\omega (A ^\flat ),\zeta_{e}\big\rangle, \quad \text{for $i,k \in N(j)$, $k \notin N(i)$},
	\end{aligned}
	\end{equation}  
	for the constant $K _i ^e $ and the vorticity $ \omega(A ^\flat )$, defined as, respectively
	\begin{equation}\label{equ_K}
	K^e_k:=\frac{|\zeta_{e}\cap T_k|}{|\zeta_{e}|} \quad \text{and} \quad 
	\big\langle\omega ({A^\flat}),\zeta_e\big\rangle:=  \sum_{ h _{mn}\in \partial \zeta _e } A^\flat _{mn} \, .
	\end{equation}
	Here, $e$ denotes the node common to triangles $T_i, T_j, T_k$ and $|\zeta_e\cap T_k|$ is the area of the intersection of 
	$T_k$ and $ \zeta _e$, where the latter denotes the dual cell to $e$, see Fig.~\ref{fig_Notation}.
	The discrete vorticity, see the second equation in \eqref{equ_K}, is calculated by taking the sum over the dual edges in the boundary $ \partial \zeta _e $ counterclockwise around node $e$. The definition of $ A ^\flat $ in \eqref{def_flat} leads to a skew-symmetric matrix, hence $ A ^\flat \in \Omega _d ^1 (\mathbb{M}  )$.				
	
	As shown in \cite{baue17a}, for the discrete Lagrangian \eqref{discrete_Lagrangian} on a simplicial grid, the discrete variational equations  \eqref{discrete_EP_compressible} lead to the following semi-discrete equations
	\begin{equation}\label{simplicial_scheme}
	\left\{
	\begin{array}{l}
	\displaystyle \vspace{0.2cm}\overline{D}_{ij} \frac{d}{dt} A^\flat_{ij} +  \omega _+ \left(K_i^+\overline{D}_{ji_+} A_{ii_+}+ K_j^+\overline{D}_{ij_+}A_{jj_+}\right) -\omega _- \left(K_i^-\overline{D}_{ji_-} A_{ii_-}
	+ K_j^-\overline{D}_{ij_-}A_{jj_-}\right) \\
	\displaystyle \vspace{0.2cm} \qquad \qquad + \overline{D}_{ij} \frac{1}{2} \left( A_{ii_-} ^\flat A_{ii_-} +  A_{ii_+} ^\flat A_{ii_+}  + A_{ij} ^\flat A_{ij} - A_{ji} ^\flat A_{ji} 
	-  A_{jj_-} ^\flat A_{jj_-}- A_{jj_+} ^\flat A_{jj_+} \right) \\
	\displaystyle \vspace{0.2cm}\qquad \qquad + \overline{D}_{ij} \left(  \frac{\partial \epsilon }{\partial D_i}- \frac{\partial \epsilon }{\partial D_j}  \right) =0\\
	\displaystyle\frac{d}{dt} D_i = A_{ii_-}D_{i_-}+ A_{ii_+}D_{i_+} + A_{ij}D_{j}  - A_{ii}D_i,
	\end{array} \right.
	\end{equation}
	in which $\overline{D}_{ij}:=  \frac{1}{2} (D_i+D_j)$ denotes the average of the cell values and 
	\linebreak $ \omega _{\pm}:=\sum_{ h_{mn}\in \partial \zeta _{\pm} } (A ^\flat  _{mn}+  R^\flat _{mn})$ is the discrete absolute vorticity at the nodes $\pm$ at endpoints of the edge between cells $i$ and $j$, see Fig.~\ref{fig_Notation}. This is the discrete version of the RSW equations \eqref{one_form_RSW}.
	
	\subsection{Semi-discrete RSW scheme in terms of the discrete velocity field}
	
	%
	
	With a suitable choice of structure preserving time discretization, Equations~\eqref{simplicial_scheme} provide a set of fully discrete equations which can be implemented as they stand. However, as it is more familiar in the GFD community to work with velocity quantities, we proceed in rewriting the equations correspondingly.
	
	From the original definitions~\eqref{approx_A_ij} and the flat operator~\eqref{def_flat}, we 
	find the following relation between one-forms $A^\flat$, Lie algebra elements $A$ and the normal 
	velocity degrees of freedom $V_{ij}$ on the triangle edges' midpoints:
	\begin{equation}\label{equ_explicit_matrixA}
	\begin{aligned}  A_{ij} &= - \frac{1}{2 \Omega _{ii} } f_{ij} V _{ij}, \quad  \text{for all} \ j \in N(i), \ j\neq i, \\
	A_{ii} &= - A_{ij} -  A_{ii_-} - A_{ii_+}   = \sum_{k=j,i_-,i_+}\frac{1}{2 \Omega _{ii} } f_{ik} V _{ik} =:\frac{1}{2} \operatorname{div} (V)_{i} .
	\end{aligned}
	\end{equation}
	while $A ^\flat _{ij}=- h _{ij} V _{ij} $ for $i \in N(j)$, see Fig.~\ref{fig_Notation} for the index notations. 
	Note that $\operatorname{div}(V)$ coincides with the natural finite volume divergence operator on a triangular mesh (see e.g. \cite{BauerPhD:2013}).

	\paragraph{Momentum equation.}
	The semi-discrete momentum equation in matrix--vector notation reads
	\begin{equation}\label{equ_swe_discrete_expl}
	\partial _t    V_{ij} +  \operatorname{Adv}_{\rm}(V,D)_{ij}= \operatorname{K}_{\rm}( V )_{ij} - \operatorname{G}(D)_{ij},
	\end{equation} 
	where we define
	\begin{linenomath}
		\begin{align*} 
		& \operatorname{Adv}_{\rm }(V,D)_{ij} : =  \\ 
		& \;\;  -  \frac{1}{\overline{D}_{ij}  h_{ij} } \Big( \frac{1}{ |\zeta_- | } 
		\sum_{ h _{mn}\in \partial \zeta _-} 
		h_{mn} { (V_{mn}  + \bar R_{mn}) } \Big) \left(   \frac{|\zeta_{-} \cap T_i|}{2 \Omega_{ii}  } \overline{D}_{ji_-}  f_{ii_-}   V_{ii_-} 
		+  \frac{|\zeta_{-}  \cap T_j|}{2 \Omega_{jj}  } \overline{D}_{ij_-}  f_{jj_-}   V_{jj_-} \right) \\
		&\;\; +  \frac{1}{\overline{D}_{ij}  h_{ij} } \Big( \frac{1}{ |\zeta_ + | }\sum_{ h _{mn}\in \partial \zeta _+} 
		h_{mn} { (V_{mn}  + \bar R_{mn}) } \Big) \left(  \frac{|\zeta_{+}  \cap T_i|}{2 \Omega_{ii} }\overline{D}_{ji_+} f_{ii_+}   V_{ii_+}
		+ \frac{|\zeta_{+}  \cap T_j|}{2 \Omega_{jj}   } \overline{D}_{ij_+}   f_{jj_+}   V_{jj_+} \right), \\
		& \operatorname{K}_{\rm }( V) _{ij}:= -  \frac{1}{2 h_{ij} }  \Big(  \frac{h_{jj_-} f_{jj_-} (V_{jj_-})^2  }{2\Omega_{jj}}  
		+ \frac{h_{jj_+} f_{jj_+} (V_{jj_+})^2  }{2\Omega_{jj}} 
		+ \frac{h_{ij  } f_{ij  } (V_{ji  })^2  }{2\Omega_{jj}} \\
		&  \qquad \qquad  \qquad      - \frac{h_{ii_-} f_{ii_-} (V_{ii_-})^2  }{2\Omega_{ii}}
		- \frac{h_{ii_+} f_{ii_+} (V_{ii_+})^2  }{2\Omega_{ii}} 
		- \frac{h_{ij  } f_{ij  } (V_{ij  })^2  }{2\Omega_{ii}} \Big)   , \\
		&\operatorname{G}(D)_{ij} :=  \frac{g}{ h_{ij} }  \Big( D_j + B_j - (D_i +  B_i)   \Big)   ,
		\end{align*}
	\end{linenomath}
	
	
	for values $\bar R_{mn}$ related to $ R_{mn}$ by $ R_{ij} = - \frac{1}{2 \Omega _{ii} } f_{ij} \bar R_{ij}$,
	analogously to the relation between $V_{mn}$ and $A_{mn}$. 
	We again use the definition $\overline{D}_{ij} = \frac{D_i + D_j }{2}$.	
	We define as Coriolis parameter
	\begin{equation}
	f|_{ \zeta _\pm}:= \frac{1}{ |\zeta_\pm | } \sum_{  h _{mn}\in \partial \zeta _\pm} h_{mn}  \bar R_{mn}.
	\end{equation}

	\paragraph{Continuity equation.}
	
	The semi-discrete continuity equations in matrix--vector form is given by 				
	\begin{equation}\label{equ_disc_cont}
	\begin{array}{l}
	\partial _t D _i  + \underbrace{\frac{1}{\Omega _{ii}} f _{ij}  V _{ij} \overline{D}_{ij}   + \frac{1}{\Omega _{ii}} f _{ii_-}  V _{ii_-} \overline{D}_{ii_-}
		+ \frac{1}{\Omega _{ii}} f _{ii_+}  V _{ii_+} \overline{D}_{ii_+}}_{:= {\rm div} (V, {D})_i}  = 0 .\\ 
	\end{array}
	\end{equation}
	Hence, the spatial variational discretization process leads to a standard finite volume representation of the 
	divergence operator and hence of the continuity equation.

	\paragraph{Time discretizations.} Since the spatial discretization has been realized by variational principles in a structure-preserving way, a temporal variational discretization can be implemented by following the discrete (in time) Euler--Poincar\'e--d'Alembert approach, analogously to what has been done in \cite{GaMuPaMaDe2011} and \cite{DeGaGBZe2014}, to which we refer for a detailed treatment. This variational approach is based on the introduction of a local approximant to the exponential map of the Lie group, see \cite{BRMa2009}, chosen here as the Cayley transform. As explained in \cite{GaMuPaMaDe2011,DeGaGBZe2014}, by dropping cubic terms, this results in a Crank--Nicolson-type time update for the momentum equation \eqref{discrete_EP_compressible} \changed{(which, as such, is second order in time)} and an update equation based on the Cayley transform for the advection equation \eqref{continuity}. Following \cite{baue17a}, we will use below the Crank--Nicolson-type time update
	directly on the momentum equation as reformulated in \eqref{equ_swe_discrete_expl}. This considerably simplifies the solution procedure without altering the behavior of the scheme. 				
	
	Alternatively, we apply also a standard time integrator using a Crank--Nicolson-type time update for the continuity equations instead of the Cayley transform and compare both time stepping schemes with each other while keeping the spatial variational discretization unmodified.				
	\medskip
	
	\noindent\underline{1.) Fluid depth equation update by Cayley transform:}
	This time integrator consists of two steps. 
	We first compute the update equation for the fluid depth $D$, which is based on the Cayley transform $\tau$. 
	This update equation is then given by $D^{t+1}= \tau (\Delta t A^t)D^t$
	for the time $t$ and a time step size $\Delta t$. In particular, $\tau$ can be represented as
	\begin{equation}\label{equ_D} 
	\big(I - \frac{1}{2} \Delta t A^t\big) D^{t+1} = \big(I + \frac{1}{2} \Delta t A^t\big) D^{t},
	\end{equation}
	with $I$ the identity matrix (cf. \cite{DeGaGBZe2014} for more details). 
	Note that $A$ can be expressed in terms of~$V$ using \eqref{equ_explicit_matrixA}.
	In a second step, we solve the momentum equation, given by an implicit nonlinear equation (step 2),
	according to the fixed-point iteration:
	\begin{enumerate}
		\item Start loop over $k=0$ with initial guess at $t$: $ V^{*}_{k=0} = V^{t}$;
		\item Calculate updated velocity $V_{k+1}^{*}$ from the explicit equation:
		\[
		\frac{V^{*}_{k+1} -V^{t}}{\Delta t}  =-\frac{\operatorname{Adv}_{\rm }(V^*_k,D^{t+1})+\operatorname{Adv}_{\rm }(V^{t  },D^{t  })}{2} +\frac{\operatorname{K}_{\rm }(V^*_k )+\operatorname{K}_{\rm }(V^{t  } )}{2}-\operatorname{G}(D^{t+1});
		\]
		\item Stop loop over $k$ if $||V^{*}_{k+1} - V^{*}_k|| < \epsilon$ for a small positive $\epsilon$, take $V^{t+1} = V^{*}_{k+1}$.
	\end{enumerate}
	
	
	For more details, we refer the reader to \cite{baue17a}.
	Note that for this time integration scheme, we do not discretize the continuity equation~\eqref{equ_disc_cont} directly,
	but use the discretization of $\tau$ in \eqref{equ_D}.
	\ \\
	
	\noindent
	\underline{2.) Fluid depth equation update by Crank--Nicolson:}
	Here, we use a two-step time integration scheme 
	to solve the system of fully discretized nonlinear momentum and continuity equations: 
	\begin{linenomath}
		\begin{align}
		\frac{V^{t+1}_{ij} -V^{t}_{ij} }{\Delta t} &  =-\frac{\operatorname{Adv}_{\rm }(V^{t+1},D^{t+1})_{ij}+\operatorname{Adv}_{\rm }(V^{t  },D^{t  })_{ij}}{2}  
		+\frac{\operatorname{K}_{\rm }(V^{t+1} )_{ij} +\operatorname{K}_{\rm }(V^{t  } )_{ij}}{2}  										- \operatorname{G}(D^{t+1})_{ij} \, ,\label{disc_update_mom} \\
		\frac{D^{t+1}_i    - D^{t}_i }{\Delta t}   &  = -\frac{\operatorname{ div} (V^{t+1}, {D}^{t+1})_i+\operatorname{ div} (V^{t  }, {D}^{t  })_i }{2}\,.\label{disc_update_cont}
		\end{align}
	\end{linenomath}
	We solve this system of nonlinear equations by fixed-point iteration 
	for all edges $ij$ and cells~$i$. To enhance readability, we skip the corresponding subindices in the following.
	The solution algorithm reads:
	\begin{enumerate}
		\item Start loop over $k=0$ with initial guess at $t$: $ V^{*}_{k=0} = V^{t}$ and $D^*_{k=0} = D^t$;
		\item Calculate updated water depth (density) $D_{k+1}^{*}$ from the explicit equation:
		\[
		\frac{D^*_{k+1}    - D^{t} }{\Delta t}    = -\frac{\operatorname{ div} (V^{*}_k,D^{*}_k)+\operatorname{ div} (V^{t  }, {D}^{t  }) }{2}
		\]
		\item Calculate updated velocity $V_{k+1}^{*}$ from the explicit equation:
		\[
		\frac{V^{*}_{k+1} -V^{t}}{\Delta t}  =-\frac{\operatorname{Adv}_{\rm }(V^*_k,D^*_{k+1})+\operatorname{Adv}_{\rm }(V^{t  },D^{t  })}{2} +\frac{\operatorname{K}_{\rm }(V^*_k )+\operatorname{K}_{\rm }(V^{t  } )}{2}- \operatorname{G}(D^{*}_{k+1});
		\]
		\item Stop loop over $k$ if $||V^{*}_{k+1} - V^{*}_k|| + ||D^{*}_{k+1} - D^{*}_k||< \epsilon$ for a small positive $\epsilon$.
	\end{enumerate}
	Note that in case of convergence, i.e. $V^{*}_{k+1} \rightarrow V^{t+1}$ and $D^{*}_{k+1} \rightarrow D^{t+1}$,
	this algorithm solves equations~\eqref{disc_update_mom} and \eqref{disc_update_cont}.

	\section{Numerical analysis of the differential operators}\label{sec_notations}
	
	We present a convergence study of the gradient, divergence and curl operators on the icosahedral meshes. These operators are used directly and indirectly in our scheme. The study will be done in $x, y, z$ coordinates on $\mathbb{R}^3$, such that it is independent of local coordinates of the hypersurface on which we solve the equations. We consider the Euclidean metric $\langle\,,\rangle$ on $\mathbb{R}^3$ and denote by $\nabla$ the gradient relative to it.
	\\~\\
	Let $\mathbf{N_x}$ be the outward unit normal vector on the sphere at $\mathbf x\in\mathcal{S}$ 
	and $P_{\mathbf x}=I-\mathbf{N_x}\mathbf{N_x}^\top$ be the orthogonal projection onto the tangent space of $\mathcal{S}$ at $\mathbf x$, where $I$ is the $3\times 3$ identity matrix. In the following, we assume that $g\colon \mathcal{S}\rightarrow\mathbb{R}$ is a given real-valued function on $\mathcal S$, with $\tilde g\colon\mathbb{R}^3\rightarrow\mathbb{R}$ being an arbitrary extension of $g$ to $\mathbb{R}^3$, and $\mathbf{u}$ is a given vector field on $\mathcal{S}$, i.e.\ $\mathbf{u}\colon\mathcal{S}\rightarrow T\mathcal{S}$, with $\tilde{\mathbf{u}}\colon\mathbb{R}^3\rightarrow\mathbb{R}^3$ being an arbitrary extension of $\mathbf{u}$ to $\mathbb{R}^3$. From the general expression of the covariant derivative induced on hypersurfaces, see e.g., \cite{Lee1997}, the gradient, divergence and curl operators on $\mathcal S$ relative to $\gamma$ can be written as
	\begin{linenomath}
		\begin{align*}
		\operatorname{grad} g(\mathbf x)  &:= \nabla \tilde g(\mathbf x) -(\nabla \tilde g(\mathbf x)\cdot \mathbf{N_x})\mathbf{N_x}=(I-\mathbf{N_x}\mathbf{N_x}^\top)\nabla \tilde g(\mathbf x) = P_{\mathbf x} \nabla \tilde g(\mathbf x)
		\\
		\operatorname{div} \mathbf{u}(\mathbf x)&:=  \langle P_{\mathbf x}\nabla,  \tilde{\mathbf{u}} (\mathbf x)\rangle
		\\
		\operatorname{curl} \mathbf{u} (\mathbf x) &:=\langle(P_{\mathbf x}\nabla)\times  \tilde{\mathbf{u}}), \mathbf{N_x}\rangle,                           
		\end{align*}    \end{linenomath}see also~\cite{flye09a,flye12a}.
	
	Note that in contrast to the numerical schemes derived in~\cite{flye09a,flye12a} we use the above analytical expressions only to initialize the initial conditions for the numerical test cases presented in Section~\ref{sec:NumericalSimulationsVariationalIntegrator}, and for comparing against the numerically obtained expressions for $\textup{grad}$,  $\textup{div}$ and $\textup{curl}$.

	As explained further below, in our variational scheme \eqref{equ_swe_discrete_expl}--\eqref{equ_disc_cont}, the discretization of these operators are given as follows
	\begin{linenomath}
		\begin{align}
		(\operatorname{Grad}^\text{num} g_{\mathcal{}})_{ij}&=\frac{g_{j}-g_{i}}{h_{ij}}, \label{full_gradient_1}\\
		(\Div^\text{num} \mathbf{u}_{\mathcal{}})_i&=\frac 1 {\Omega_{ii}} \sum_{\ell \in \{j,i_-,i_+\}} f_{i\ell}V_{i\ell},\label{full_div_1}\\			
		(\Curl^\text{num} 	\mathbf{u}_{\mathcal{}})_{\zeta_e}
		&=\frac{1}{|\zeta_e|}\sum_{mn\in\partial\zeta_e} h_{mn} V_{mn},\label{full_curl_1}
		\end{align}                                \end{linenomath}
	where, from \eqref{approx_A_ij} and \eqref{equ_explicit_matrixA}, we have $V_{i\ell}=\frac{1}{f_{i\ell}}\int_{D_{i\ell}}(\mathbf{u}\cdot \mathbf{n}_{i\ell})dS$ and where $g_i=\frac{1}{\Omega_{ii}}\int_{T_i}gdx$. As before, $\mathbf {n}_{ij}$ is the unit normal pointing from $T_i$ to $T_j$, $\zeta_e$ is the cell dual to a node $e$ (a hexagon or pentagon), $|\zeta_e|$ is its area. Note that in \eqref{full_div_1} the sum is over the cells adjacent to $T_i$, and in \eqref{full_curl_1} the sum is over the dual edges in the boundary $ \partial \zeta _e $ counterclockwise around node $e$.

	The gradient \eqref{full_gradient_1} appears in \eqref{equ_swe_discrete_expl} and is denoted as $G(D)$.
	Therein, in the advection term, $\operatorname{Adv}_{\rm}(V,D)$, we find the curl operator~\eqref{full_curl_1} consisting of a counterclockwise sum over the edges of a dual cell. The divergence operator~\eqref{full_div_1} with positive fluxes when pointing out of the triangles appears in the continuity equation~\eqref{equ_disc_cont}. 
	\color{black}
	
\begin{figure}[!ht]
		\centering
		\begin{subfigure}{0.6\textwidth}
			\includegraphics[width=\textwidth]{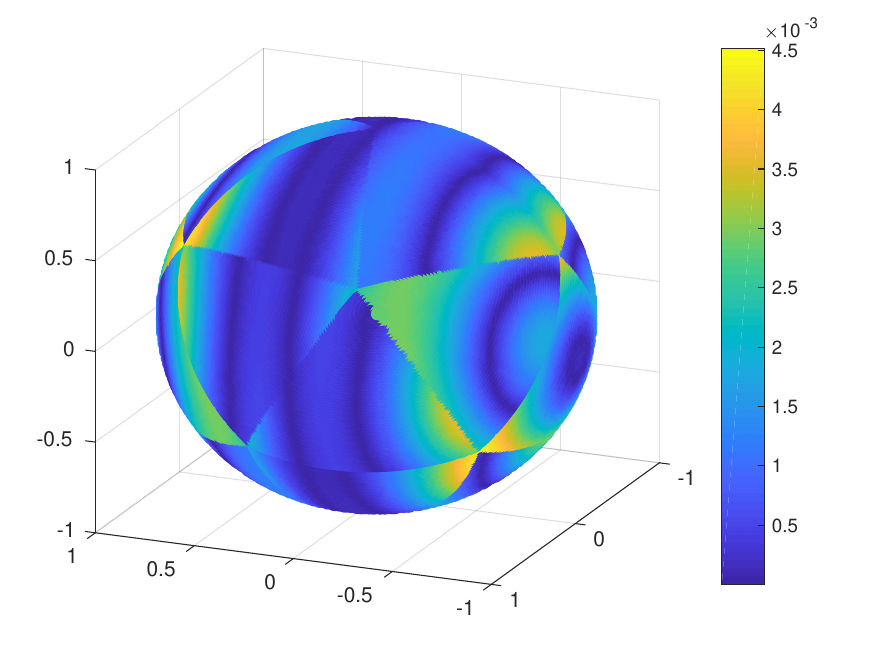}
		\end{subfigure}						
		\caption{ The grid imprint of the absolute error for the evaluation of the divergence on triangles \changed{, on a grid with
				40962 Voronoi cells (which corresponds to 81920 triangles).}}
		\label{fig:DCError}				
	\end{figure}

	For our convergence study, we use $g(\mathbf x)=\sin(x)+\sin(2y)+\sin(2z)$ for the gradient, $\mathbf{u}(\mathbf{x})=\left(x-x^3, -x^2y, -x^2z\right)^\top$ for the divergence and $\mathbf{u}(\mathbf x)=\left(z, 0,  -x\right)^\top$ for the curl.
	\begin{figure}[!ht]
		\centering
		\begin{subfigure}{0.3\textwidth}
			\includegraphics[width=\textwidth]{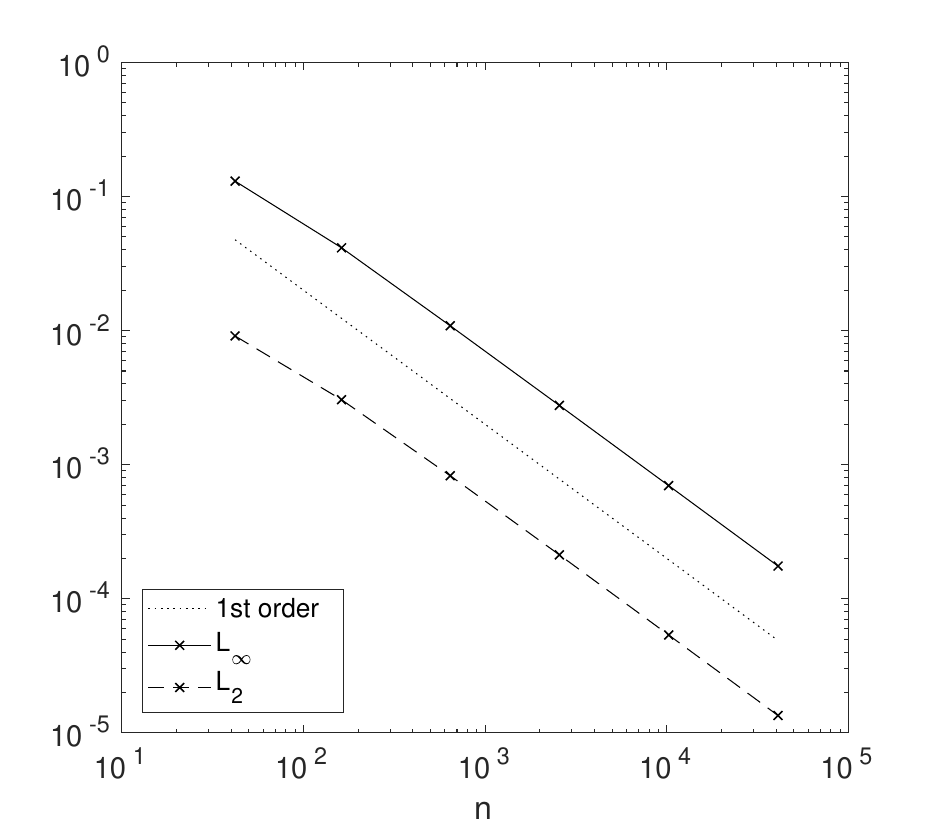}
		\end{subfigure}				
		\begin{subfigure}{0.3\textwidth}
			\includegraphics[width=\textwidth]{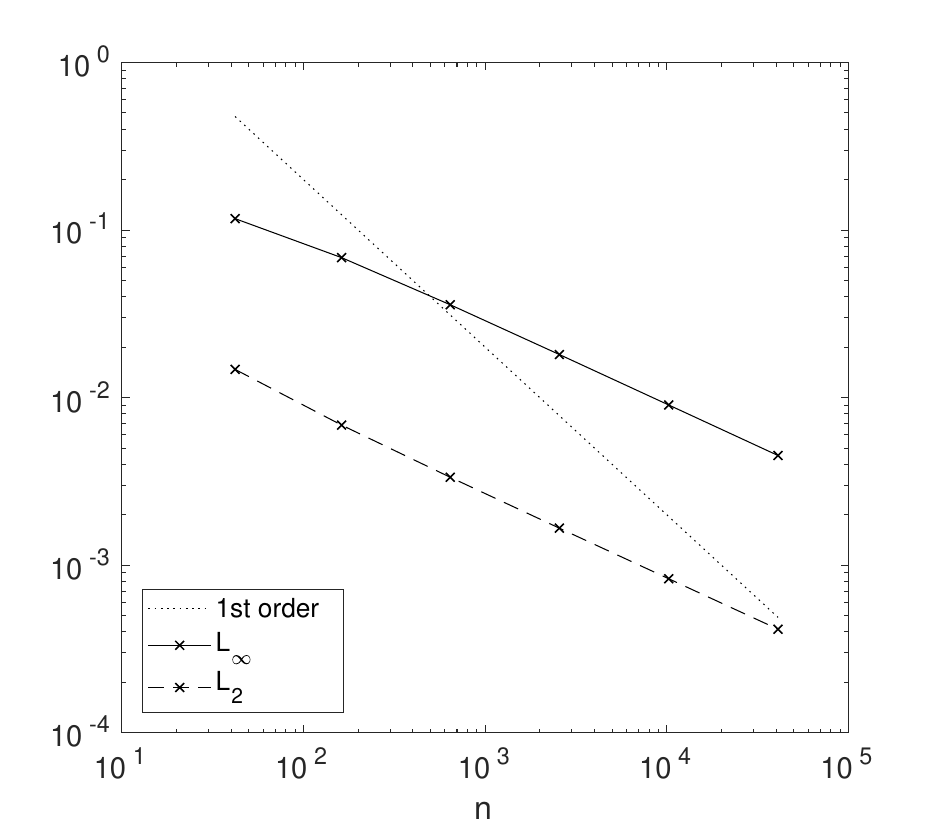}
		\end{subfigure}						
		\begin{subfigure}{0.3\textwidth}
			\includegraphics[width=\textwidth]{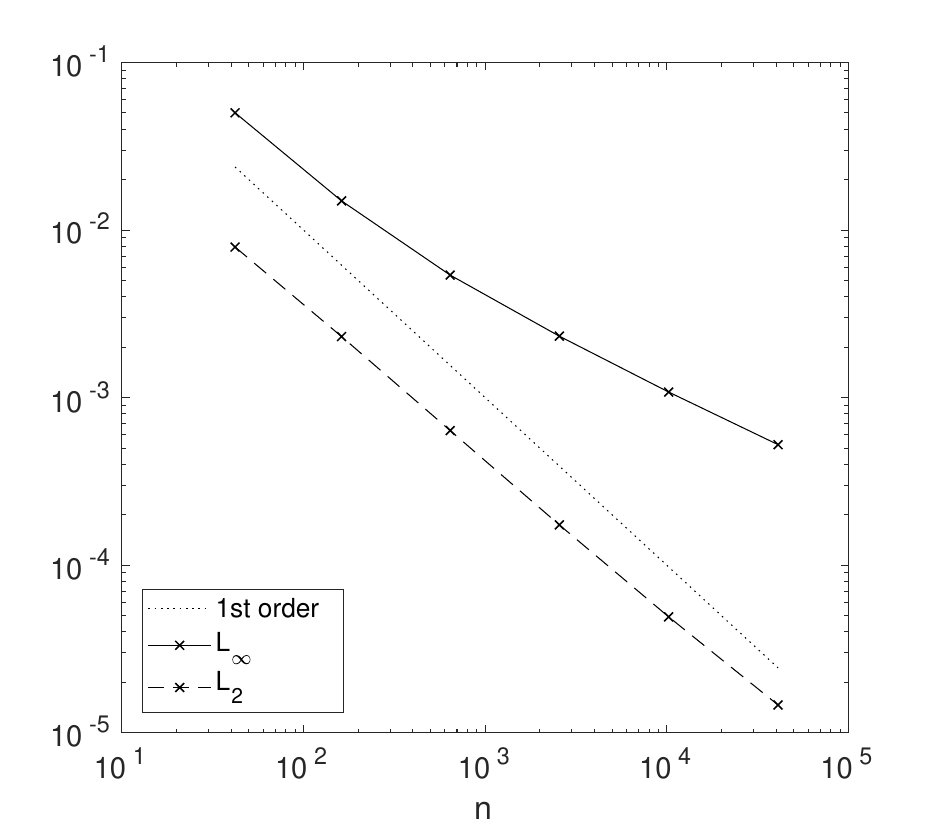}
		\end{subfigure}
		\caption{Convergence table for the discrete gradient (left), divergence (center), and curl (right).}
		\label{fig:convGDC}				
	\end{figure}\\
	In Fig.~\ref{fig:convGDC}, we see that the gradient and curl operator converge with first order. 
	The divergence operator converges too, but does not achieve first order when evaluated on the triangles, because the grid is optimized for the hexagons \cite{Heikes_Randell_II:1995}.
	In Fig.~\ref{fig:DCError}, we clearly see the grid imprint of the original icosahedron in the pointwise errors for the divergence, showing that the error is highest along the edges of the icosahedron.
	
	\section{Numerical simulations}\label{sec:NumericalSimulationsVariationalIntegrator}

	We consider four test cases; (1) the lake-at-rest solution to demonstrate that the model is well-balanced; (2) a global steady-state solution to study the convergence and energy and enstrophy loss of the model; (3) the flow over an isolated mountain; and (4) the Rossby--Haurwitz wave solution.
	
	The following constants are kept fixed for all simulations:
	\begin{linenomath}
		\[
		R= 6.37122\cdot 10^6 ~~~ [m], \qquad
		\Omega = 7.292\cdot 10^{-5} ~~~ [s^{-1}], \qquad
		g= 9.80616 ~~~ [ms^{-2}] .		
		\]                                \end{linenomath}
	The Coriolis parameter is $f=2\Omega \sin\theta$. Unless indicated otherwise, the simulations are carried out on a grid with $N=40962$ Voronoi cells (corresponding to a resolution of about 120 km) for which we chose a time step  of $\Delta t=100$ s. To estimate the numerical errors, we use the following definitions for the relative $L_\infty$-error and $L_2$-error, 	
	\begin{linenomath}						
		\begin{align*}
		\|u-u_0\|_\infty = \frac{\max_i |u(i)-u_0(i)|}{\max_i |u_0(i)|} &&
		\|u-u_0\|_2 = \frac{\sqrt{\sum_i{(\Omega_{ii}(u(i)-u_0(i)))^2}}}{\sqrt{\sum_i{\Omega_{ii} u_0(i)^2}}}\,,
		\end{align*}                     \end{linenomath}
	where $\Omega_{ii}$ is the area associated with cell $i$. Here, the function $u(i)$ is the numerical solution defined at $\mathbf x_i$ or the magnitude of the numerical solution at $\mathbf x_i$ (when used for calculating the error for the velocity) and $u_0(i)$ is the initial function at $\mathbf x_i$.
	\\~\\
	Since our scheme preserves mass and potential circulation up to machine precision, those error norms are not presented.
	In this section, we only present results obtained by using the Cayley transform time discretization, 
	which we compare in Section~\ref{sec:LongTermSimulations} with the standard time integrator.		
	\subsection{Case 1: Lake at rest}
\begin{figure}[!ht]
		\centering	
		\begin{subfigure}{0.49\textwidth}	
			\includegraphics[width=\textwidth]{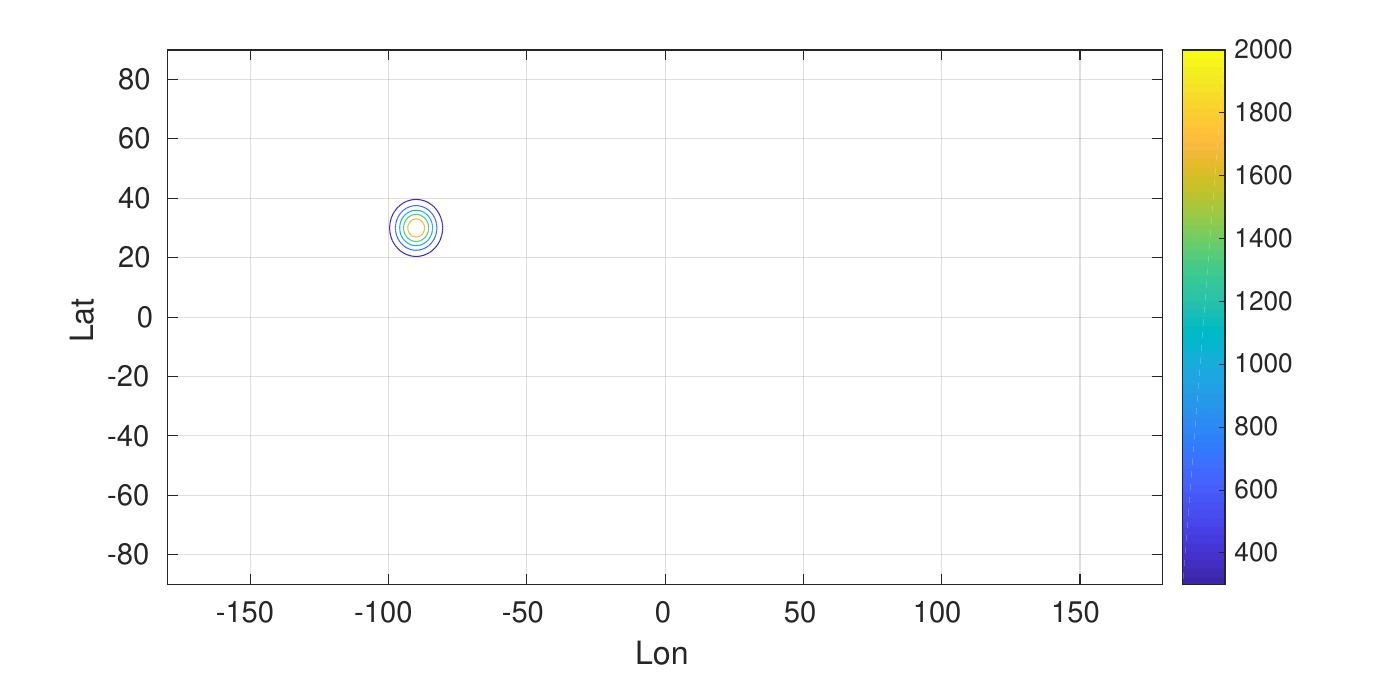}
			\caption{\footnotesize{Profile \changed{in [m]}  of the smooth bottom topography case.}}
			\label{fig:lakeMountain}
		\end{subfigure}				
		\begin{subfigure}{0.49\textwidth}	
			\includegraphics[width=\textwidth]{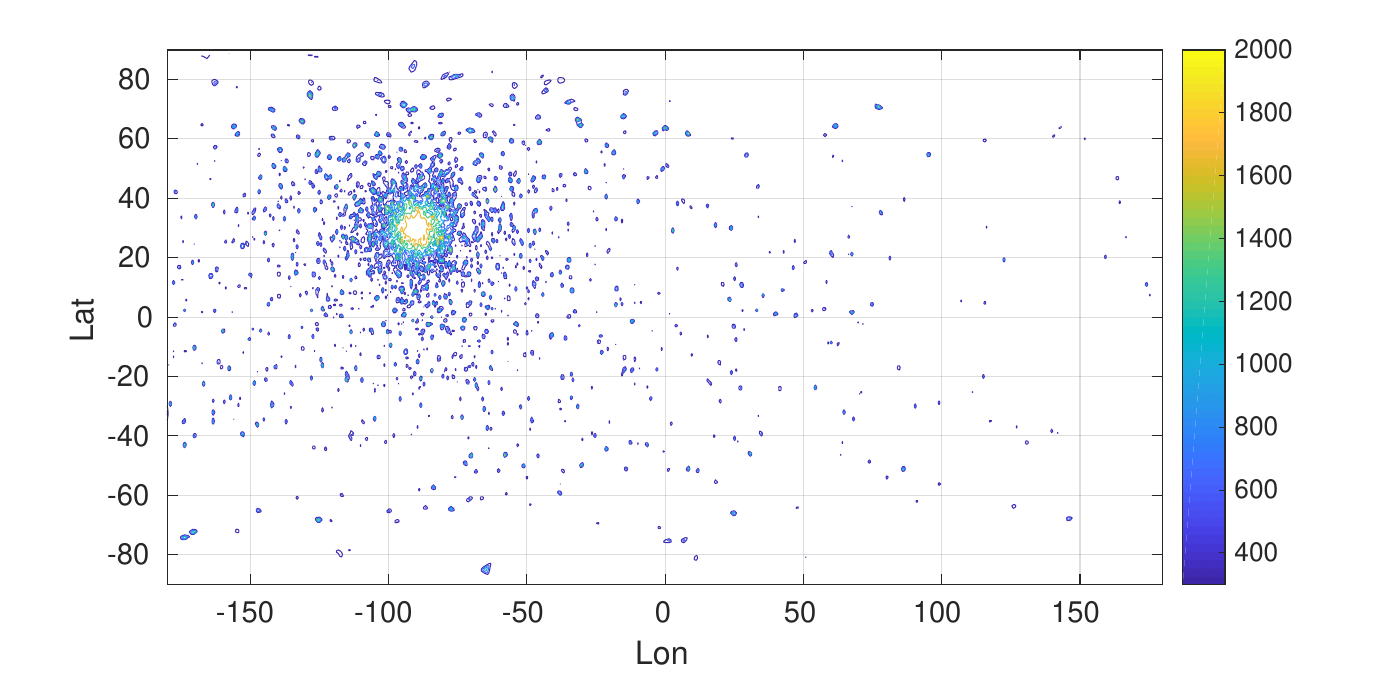}
			\caption{\footnotesize{Profile \changed{in [m]} of the noisy bottom topography case.}}
			\label{fig:lakeMountainNoise}
		\end{subfigure}	
		\caption{Topography for smooth and noisy bottom topography.}
	\end{figure}

	This test case verifies that the model is \textit{well-balanced}, that is, the exact solution $\mathbf{u}=0$, $h+B=\const$ of the RSW equations is preserved up to machine precision. Here, we choose the test case of a resting fluid over a conical shaped mountain. The initial velocity is $\mathbf{u}\equiv 0$  and the bottom topography is defined by
	\begin{linenomath}
		\[
		B(\lambda,\theta)=2000\exp\left(-(2.8\cdot9r/\pi)^2\right) \text{\quad with\quad } r^2=\min\left((\pi/9)^2,(\lambda-\lambda_c)^2+(\theta-\theta_c)^2\right),
		\]							    \end{linenomath}
	see also Fig.~\ref{fig:lakeMountain}. The total water depth is $D=5960-B$.  In addition, we carry out a second simulation where we add some white noise to the bottom topography profile, see Fig.~\ref{fig:lakeMountainNoise}, to verify that the model remains well-balanced also for a noisy bottom topography.		
	\\
	The shallow-water equations are integrated over 15 days. \changed{The initial conditions are preserved to the order of machine precision, which verifies that the model is indeed well-balanced.}  
	\subsection{Case 2: Global steady-state nonlinear zonal geostrophic flow}
	
	This test case, originally proposed in \cite{wil92}, is a geostrophically balanced flow over a flat bottom topography, i.e. $B\equiv 0~[m]$. This flow represents an exact solution to the rotating shallow-water equations. The initial conditions are: 	
	\begin{linenomath}
		\begin{align*}			
		V_{ij}=V(x_{ij}, y_{ij}, z_{ij})&=u_0(-y_{ij},x_{ij},0)^\mathsf{T}\cdot {\bf n}_{ij},   &\text{ where ~~~}&
		u_0 =\frac{2\pi R}{12\cdot 86400} ~~~[s^{-1}]
		\\				
		D_i=D(x_{T_i}, y_{T_i}, z_{T_i})&= h_0 - \frac 1 g\left(R\Omega u_0 + u_0^2/2\right)z_{T_i}^2 , &\text{ where ~~~}&
		gh_0 = 2.94\cdot 10^{4} ~~~[m^2s^{-2}].
		\end{align*}                      \end{linenomath}
	Here ${\bf n}_{ij}$ denotes the normal vector of a triangle edge, and $V_{ij}$ is the directional magnitude of the velocity normal to an edge $f_{ij}$, see Equation \eqref{equ_explicit_matrixA}.
	
	Although the nonlinear zonal geostrophic flow is a steady state solution of the RSW 
	(in which any quantity is conserved because of no time dependence), it is 
	only a stationary solution of a numerical RSW scheme up to numerical errors.
	As such, it is also interesting to monitor the time series of the numerical values of the conserved quantities of the RSW, and to verify that the energy error converges at the expected first order.
	
\begin{figure}[!ht]
		\centering
		\includegraphics[width=\textwidth]{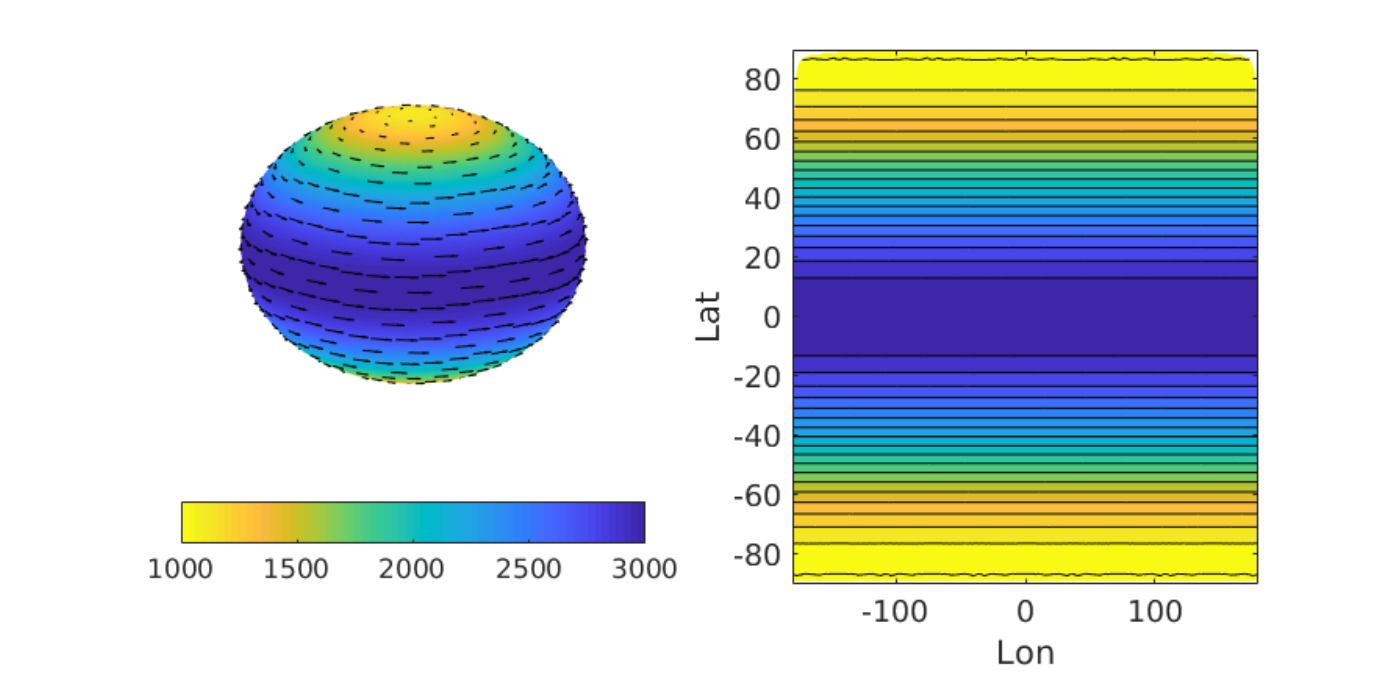}
		\caption{\footnotesize{Initial conditions for the nonlinear zonal geostrophic flow test case.}}
		\label{fig:initSteady}
	\end{figure}

\begin{figure}[!ht]
		\centering
		\includegraphics[width=\textwidth]{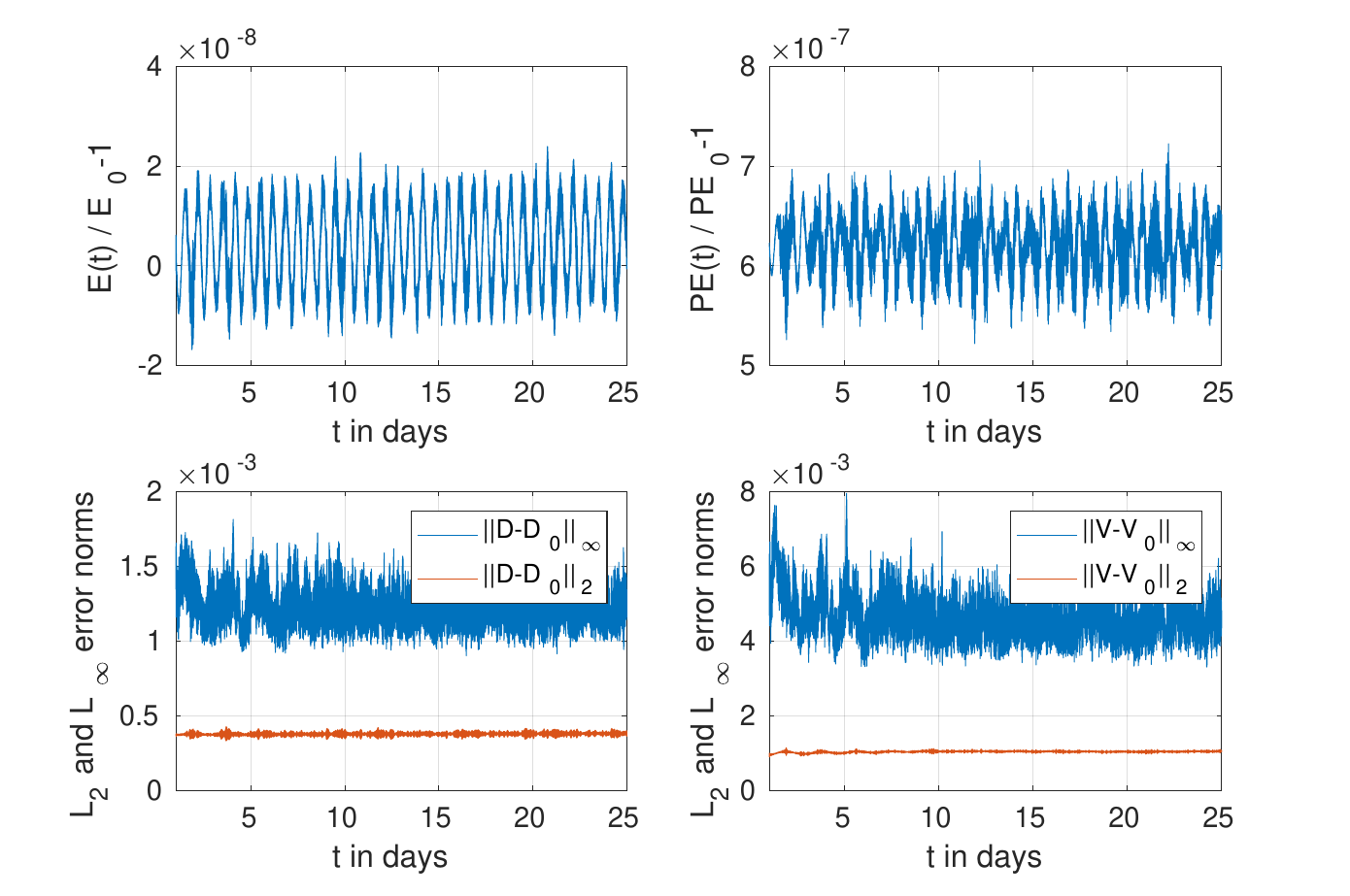}
		\caption{\footnotesize{Time series for the error norms for the conserved quantities for the nonlinear zonal geostrophic flow test case.}}
		\label{fig:errorSteady}
	\end{figure}

\begin{figure}[!ht]
		\begin{subfigure}{0.49\textwidth}
			\centering
			\includegraphics[width=\textwidth]{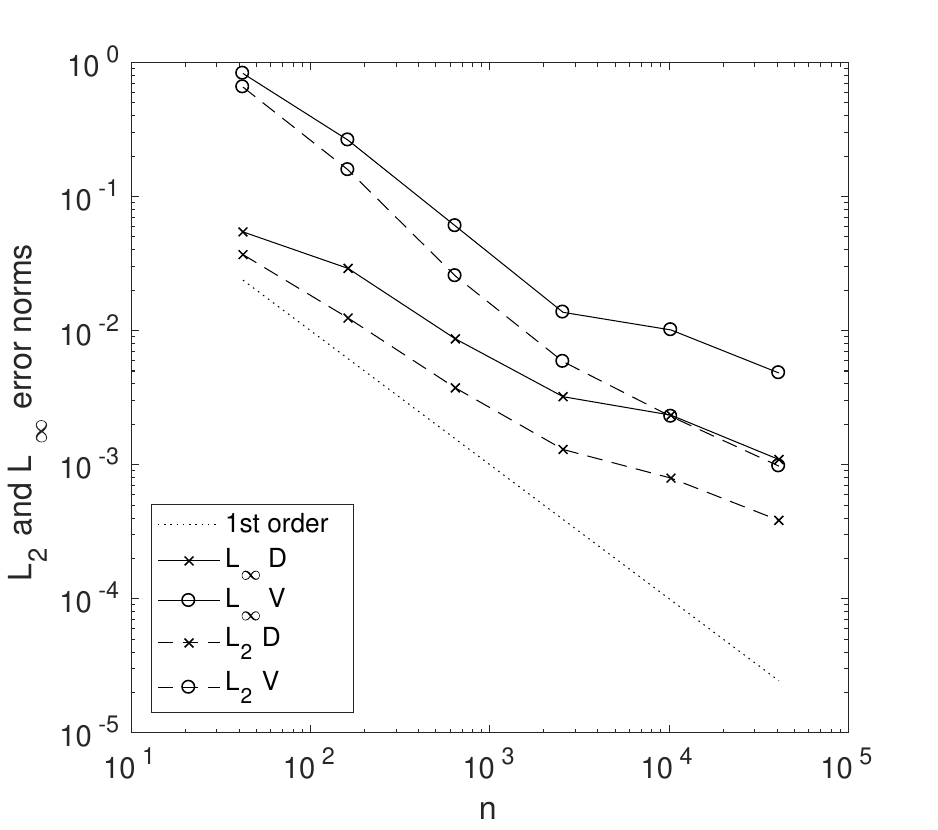}
			\caption{\footnotesize{Convergence plot of $D$ and $V$ for the nonlinear zonal geostrophic flow test case. Norms are taken at day 12.}}
			\label{fig:convSteady}
		\end{subfigure}	
		\begin{subfigure}{0.49\textwidth}
			\centering
			\includegraphics[width=\textwidth]{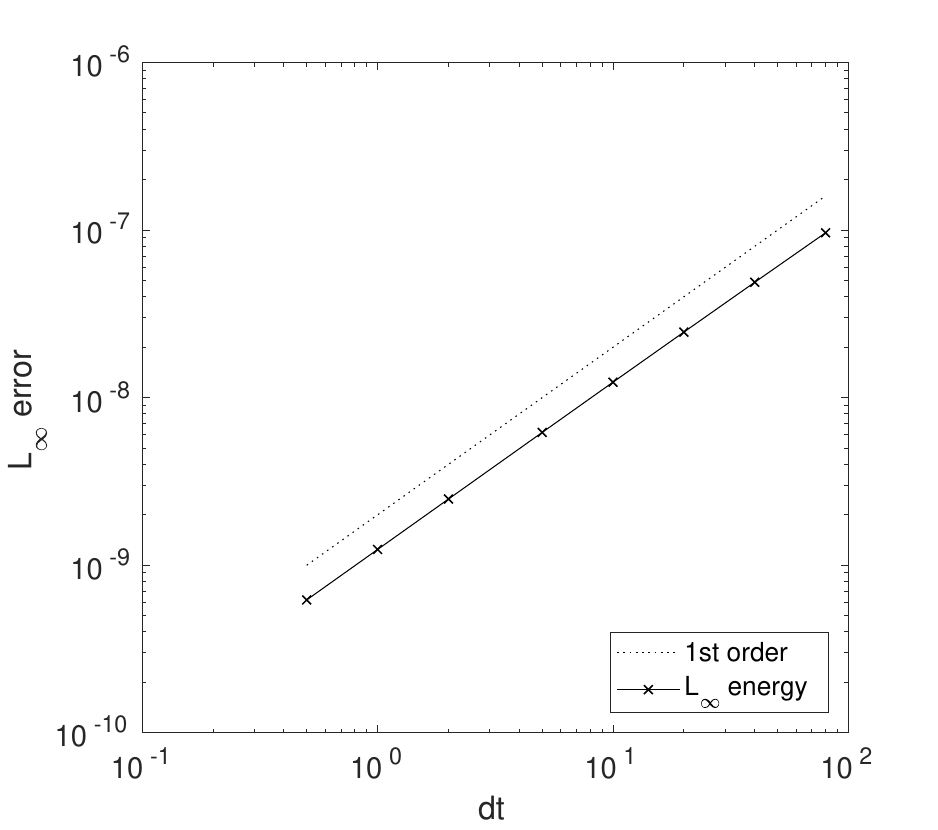}
			\caption{\footnotesize{Convergence plot of energy error for $n=10242$.}}
			\label{fig:convEnergy}
		\end{subfigure}
		\caption{Convergence plots for nonlinear zonal geostrophic flow test case}
	\end{figure}

	Fig.~\ref{fig:initSteady} shows the initial conditions. 		
	In Fig.~\ref{fig:errorSteady}, we display the time evolution of the errors for energy, potential enstrophy, height and velocity. It can be seen that there is no trend in the evolution of the error. The energy is well-preserved at the order of $10^{-8}$ and the potential enstrophy at order of $10^{-7}$. Fig.~\ref{fig:convSteady} contains the results of the spatial convergence study over 12 days with different resolutions. We integrate over 12 days, as one rotation of the fluid flow around the globe takes precisely 12 days.  It can be seen that $D$ and $V$ do not achieve first order convergence, which is natural since the divergence operator likewise does not achieve first order convergence, see Fig.~\ref{fig:DCError}. 
	Fig.~\ref{fig:convEnergy} shows the expected first order convergence of the energy error with respect to the time step.
	
	~
	
	\subsection{Case 3: Flow over an isolated mountain}				
\begin{figure}
		\centering	
		\begin{subfigure}{0.49\textwidth}	
			\includegraphics[width=\textwidth]{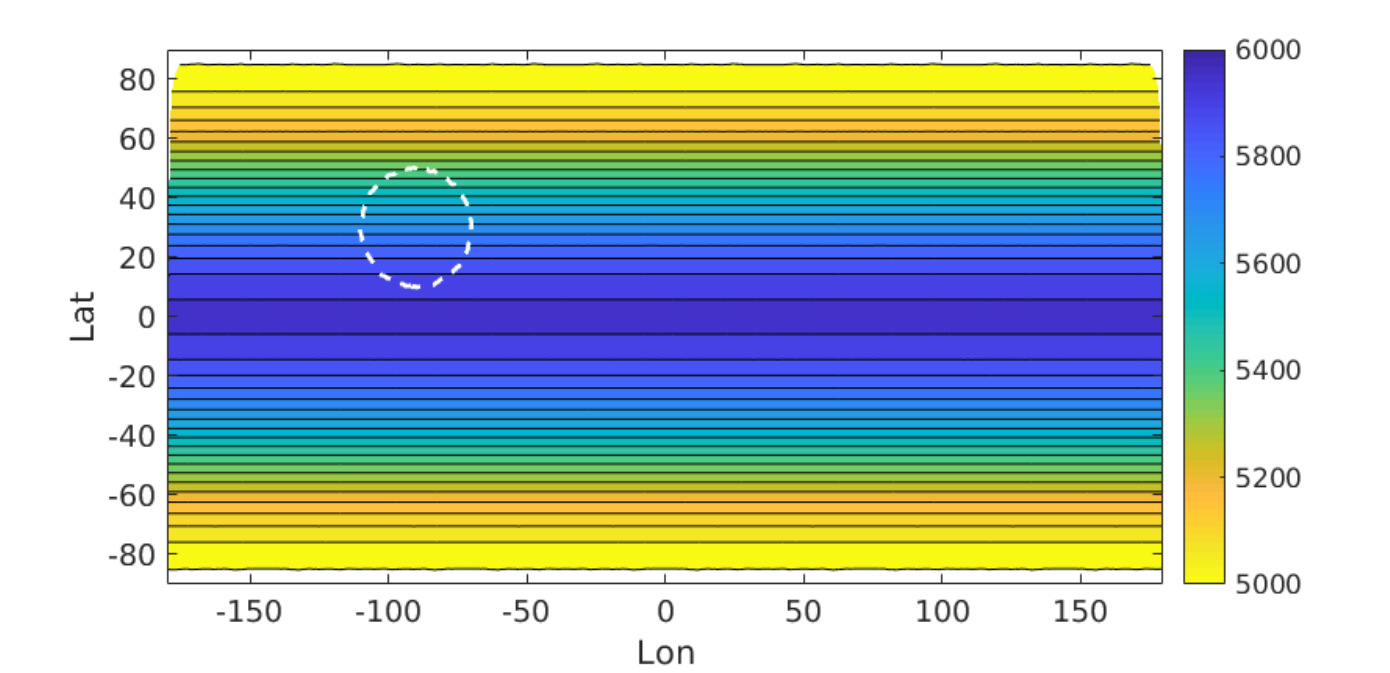}
			\caption{\footnotesize{initial}}
			\label{fig:initflow}
		\end{subfigure}					
		\begin{subfigure}{0.49\textwidth}	
			\includegraphics[width=\textwidth]{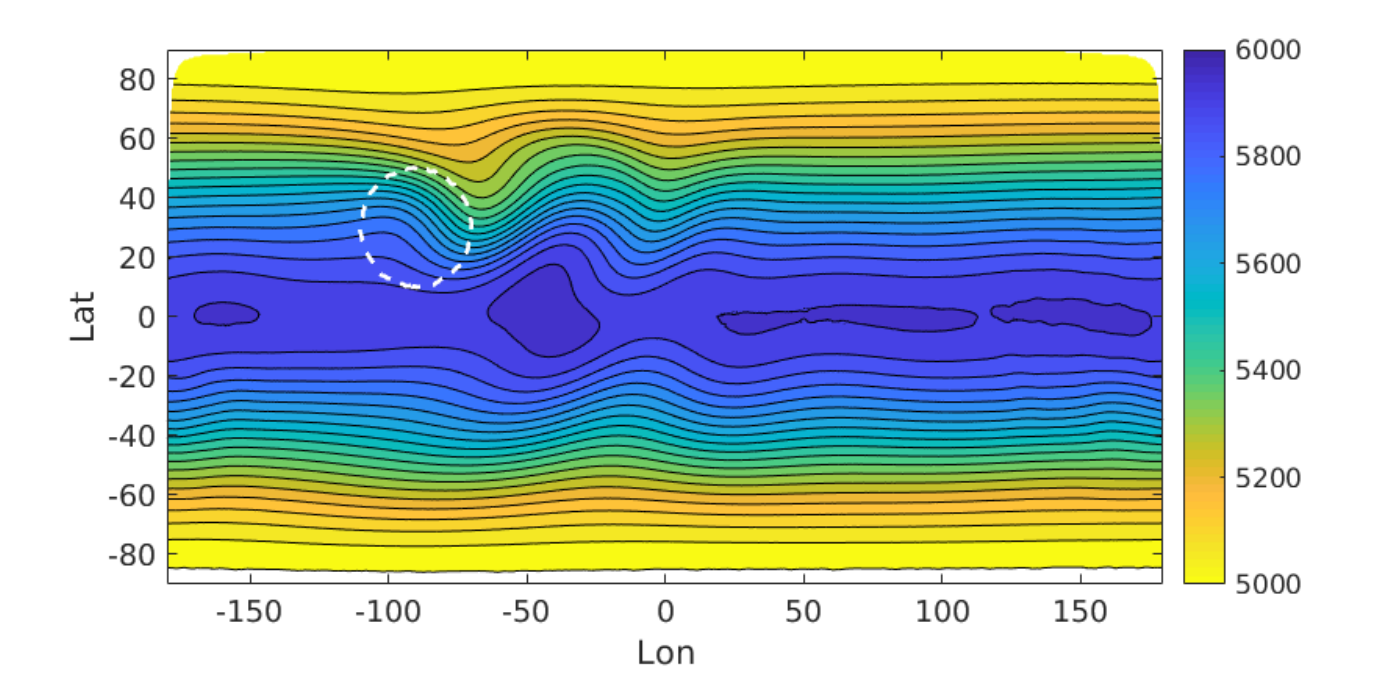}
			\caption{\footnotesize{after 5 days}}
		\end{subfigure}
		\\
		\begin{subfigure}{0.49\textwidth}	
			\includegraphics[width=\textwidth]{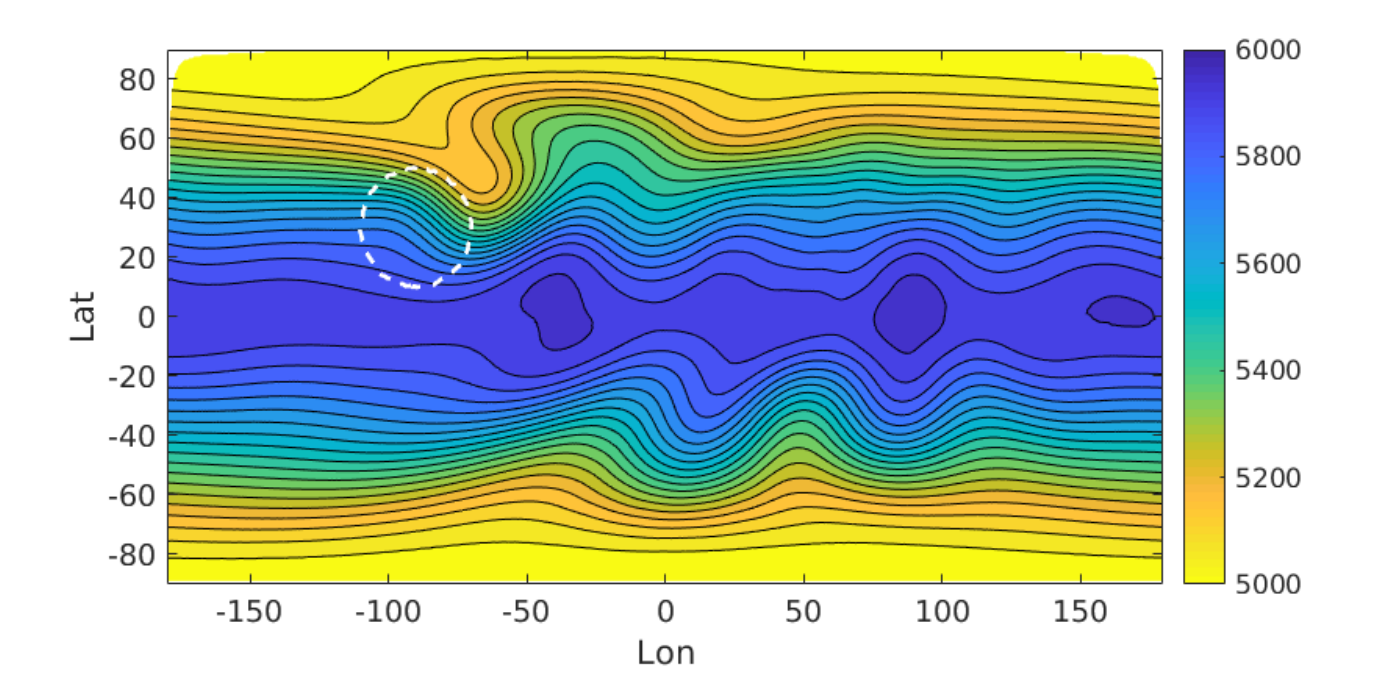}
			\caption{\footnotesize{after 10 days }}
		\end{subfigure}
		\begin{subfigure}{0.49\textwidth}		
			\includegraphics[width=\textwidth]{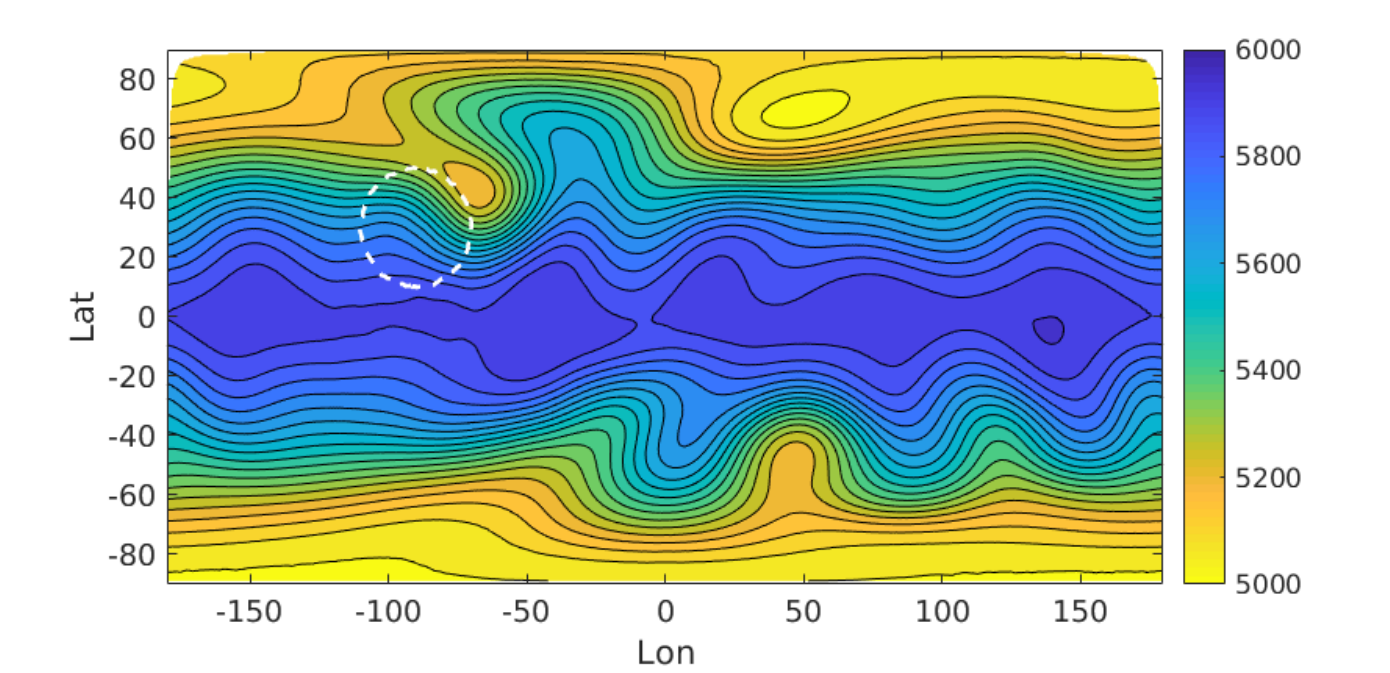}
			\caption{\footnotesize{after 15 days }}
		\end{subfigure}	
		\caption{Numerical solution for flow over an isolated mountain. Contour interval is 50m.}
		\label{fig:flow}
	\end{figure}

	Here, we consider the flow over a conically-shaped mountain which was also proposed in \cite{wil92}.
	The initial conditions, see Fig.~\ref{fig:initflow}, of this test case are the same as for Case 2 considered in the previous subsection, except that now $h_0=5960~[m]$ and $u_0=20~[m s^{-1}]$. The following bottom topography is imposed, 
	\begin{linenomath}
		\[
		B(\lambda, \theta)=2000\,(1-9r/\pi) \text{\quad with\quad } r^2=\min\big((\pi/9)^2,(\lambda-\lambda_c)^2+(\theta-\theta_c)^2\big).
		\]                                \end{linenomath}
	The mountain is centered at $\lambda_c=3\pi/2$ and $\theta_c=\pi/6$. Note that there is no analytical solution for this problem.
	
\begin{figure}[!ht]
		\includegraphics[width=\textwidth]{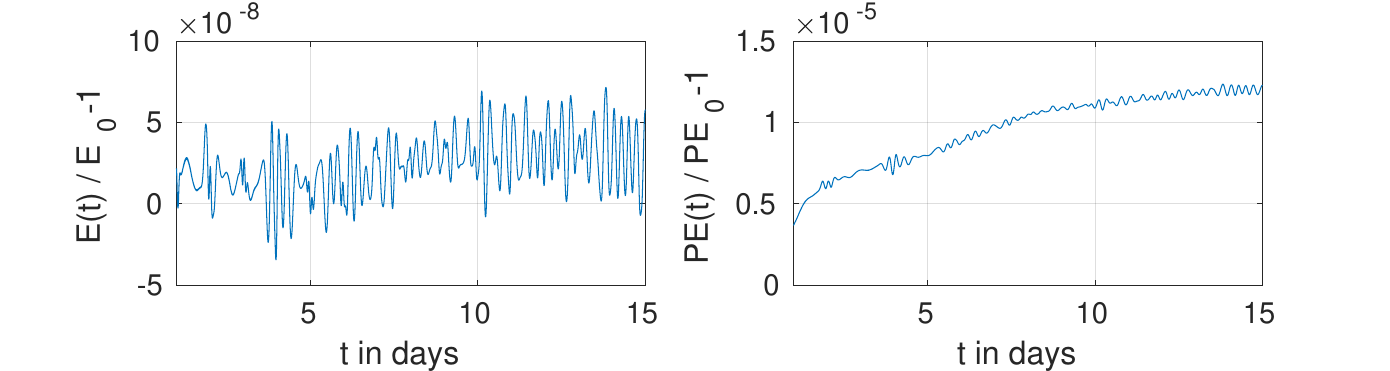}
		\caption{\footnotesize{Time series of the errors of energy and potential enstrophy for flow over an isolated mountain over the 15 days integration period.}}
		\label{fig:errorflow}
	\end{figure}

	Fig.~\ref{fig:flow} shows snapshots of the height field at times $t=0$ (a), $t=5$ days (b), $t=10$ days (c) and $t=15$ days (d), which are the times suggested in~\cite{wil92} to show the computed solutions. These results are visually similar to the results from different models, such as those given in~\cite{flye12a,ULLRICH20106104}. The time series for the errors in the energy and potential enstrophy are depicted in Figure~\ref{fig:errorflow}. 
	
	\changed{
		Note that the energy is almost as well conserved as in Case~2, whereas
		potential enstrophy conservation is three orders of magnitudes less 
		accurate (cf. Fig.~\ref{fig:errorSteady} and Fig.~\ref{fig:errorflow}).
		This is because the variational integrator preserves energy by construction, but we have no control over the conservation behavior of potential enstrophy (only PV is conserved). This, however, should not provide a problem in practice because usually enstrophy has to be dissipated anyways at small scales, cf. \cite{mcra13Ay}. Nevertheless, our schemes preserves potential enstrophy relatively well given the nonlinearity of the test case. 
	}

	\subsection{Case 4: Rossby--Haurwitz waves}
	
\begin{figure}[!ht]
		\centering
		\begin{subfigure}{0.49\textwidth}	
			\includegraphics[width=\textwidth]{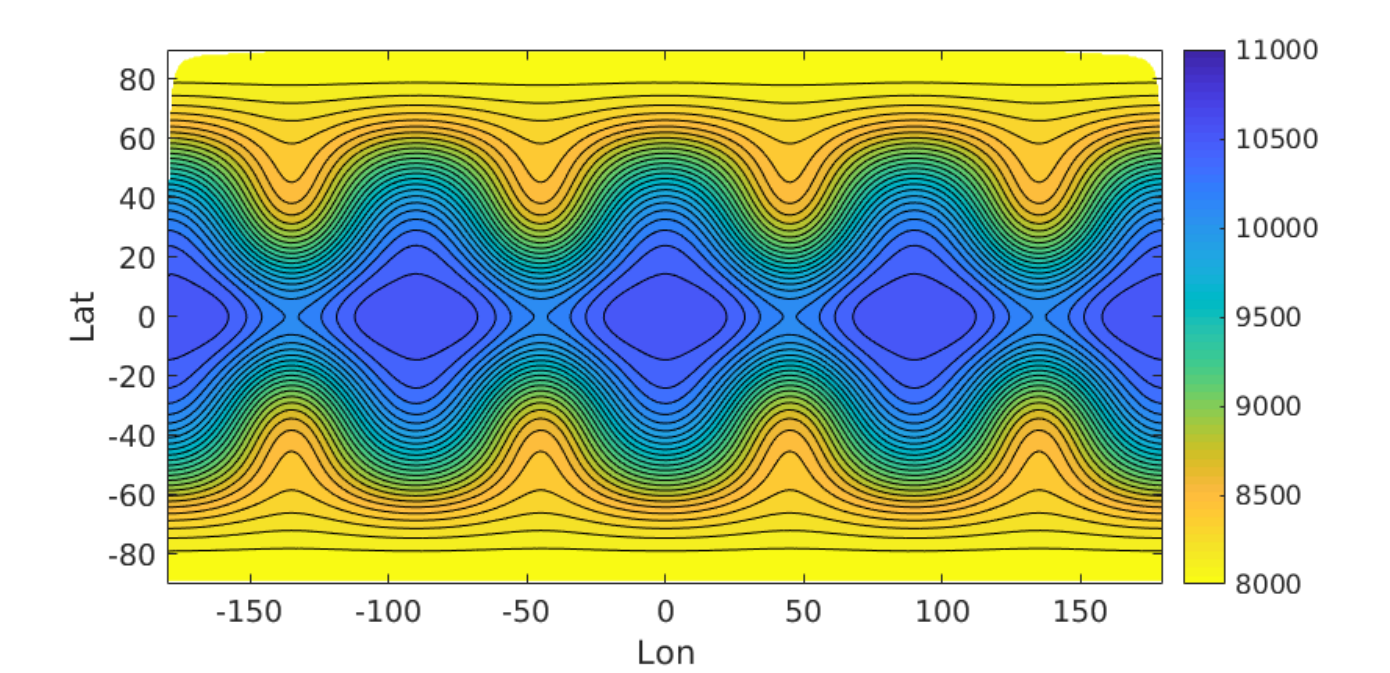}
			\caption{\footnotesize{initial wave form}}
			\label{fig:initrossby}
		\end{subfigure}		
		\begin{subfigure}{0.49\textwidth}
			\includegraphics[width=\textwidth]{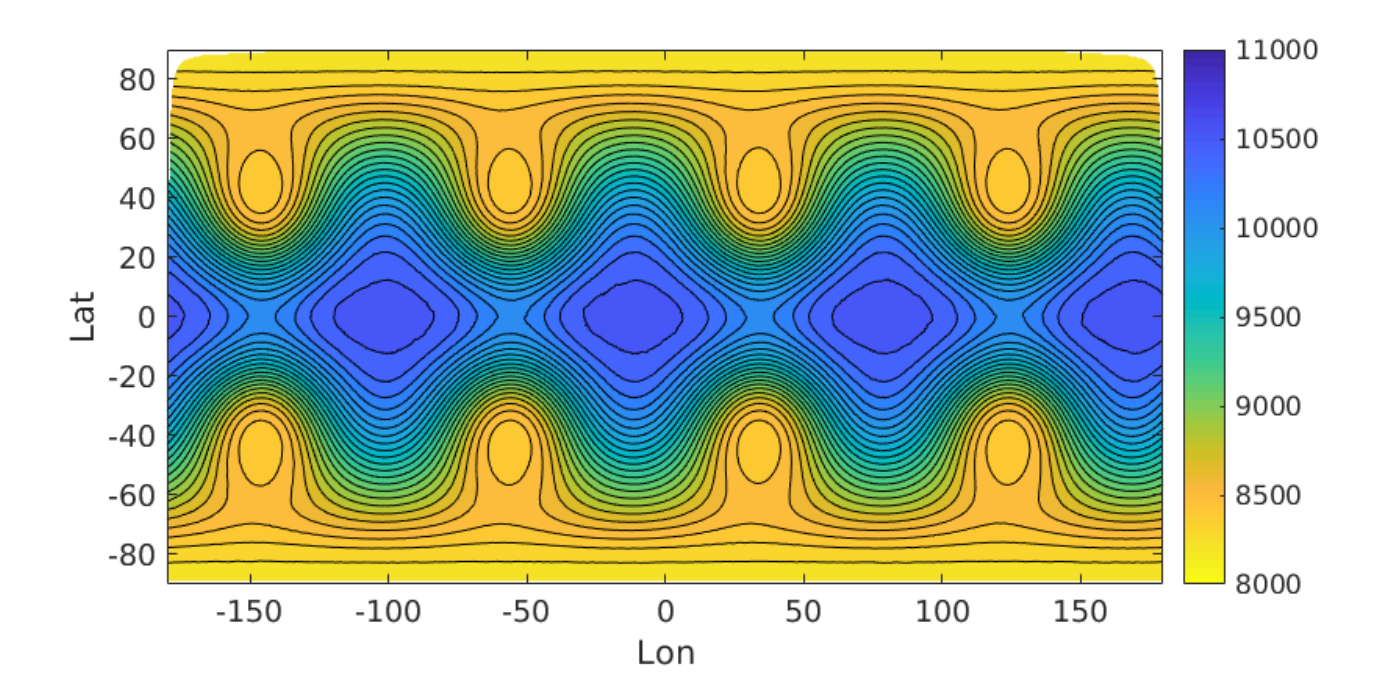}
			\caption{\footnotesize{form of the wave after 7 days} }
			\label{fig:7rossby}
		\end{subfigure}\\
		\begin{subfigure}{0.49\textwidth}
			\includegraphics[width=\textwidth]{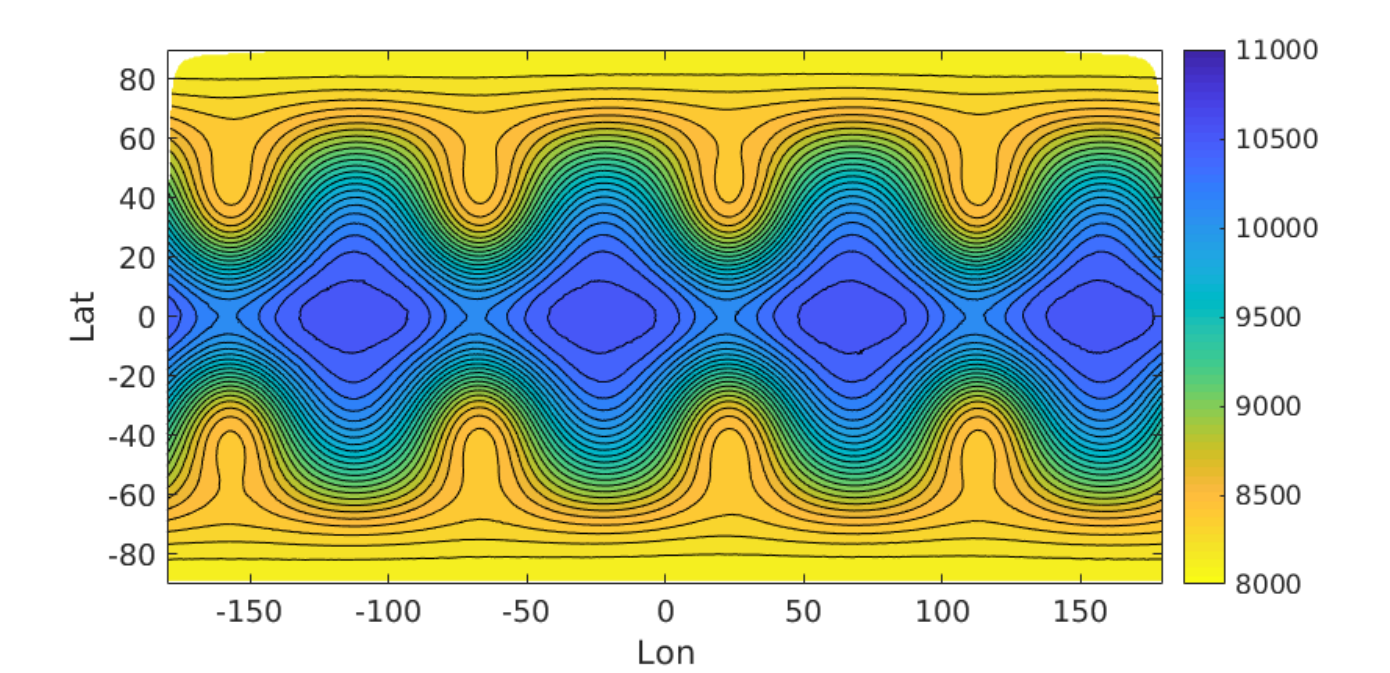}
			\caption{\footnotesize{form of the wave after 14 days}}
			\label{fig:14rossby}
		\end{subfigure}
		\begin{subfigure}{0.49\textwidth}
			\includegraphics[width=\textwidth]{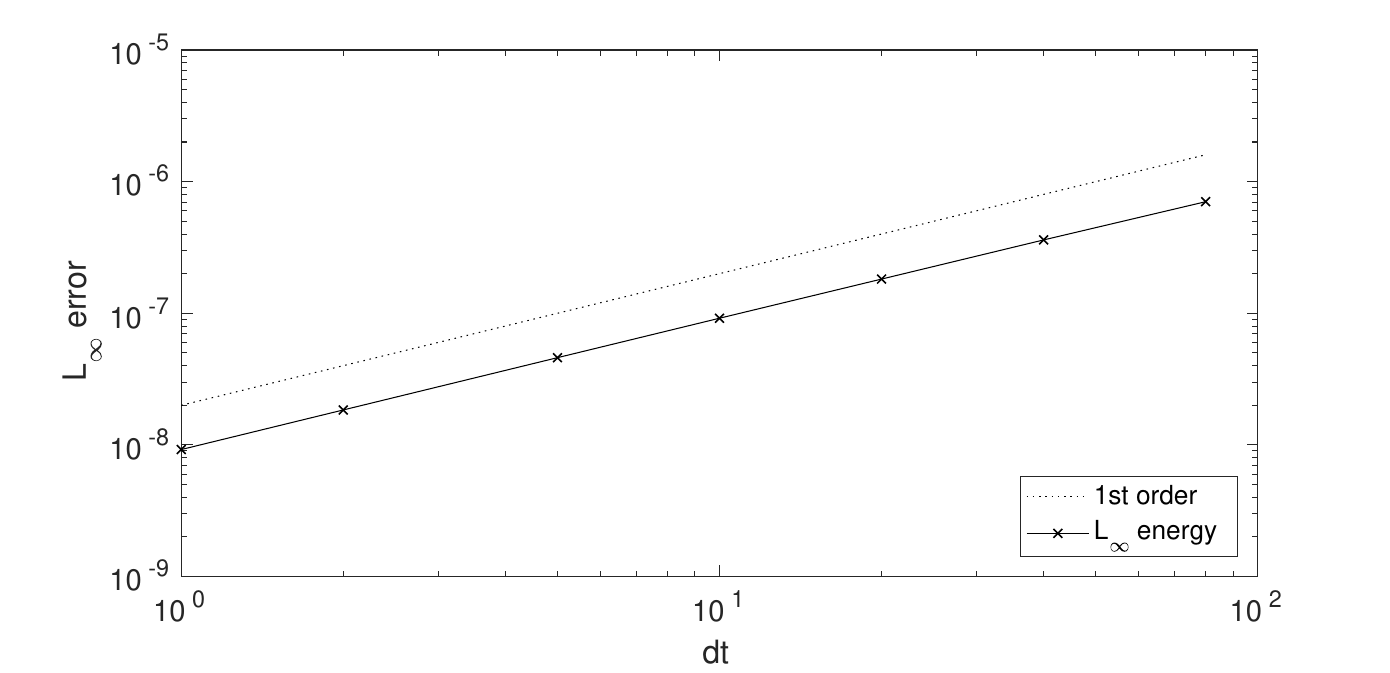}
			\caption{\footnotesize{convergence of energy error for $n=10242$}}
			\label{fig:rossbyEnergy}
		\end{subfigure}	
		\caption{ Numerical solution for Rossby-Haurwitz waves. Contour interval is 100m.}
		\label{fig:rossby}
	\end{figure}

\begin{figure}[!ht]
		\centering
		\includegraphics[width=\textwidth]{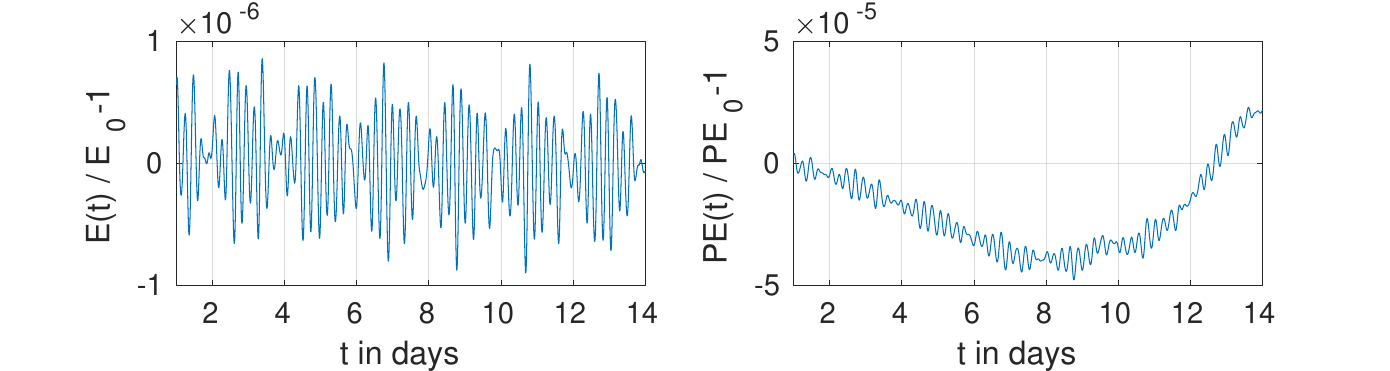}
		\caption{Time series of the errors of energy and potential enstrophy for Rossby--Haurwirz waves over the 14 days integration period.}
		\label{fig:errorrossby}
	\end{figure}

	We consider a Rossby--Haurwitz wave with wavenumber $\kappa=4$, which is proposed in \cite{wil92}. Unlike the non-divergent barotropic vorticity equation, the shallow-water equations can only approximate this solution. For comparison, snapshots after 7 and 14 days are presented. For completeness we present the initial conditions here in latitude ($\theta$) and longitude ($\lambda$),
	\begin{linenomath}
		\begin{align*}
		K=7.848\cdot 10^{-6}\quad [s^{-1}],\qquad \kappa= 4\quad (\text{wave number}),\qquad
		h_0= 8\cdot 10^{3}\quad [m].
		\end{align*}                       \end{linenomath}
	The components of the velocity vector $\mathbf{u}=(u,v)$ are
	\begin{linenomath}
		\begin{align*}
		u&=RK \cos \theta + R K \cos^{\kappa-1}\theta (\kappa \sin^2\theta -\cos^2\theta) \cos\kappa\lambda,
		\\
		v&=-RK\kappa\cos^{\kappa-1}\theta\sin\theta \sin\kappa\lambda	
		\end{align*}                     \end{linenomath}
	and the water elevation is given by
	\begin{linenomath}
		\begin{align*}
		D_{i}&=h_0+\tfrac 1 g(R^2A(\theta_{T_i})+R^2B(\theta_{T_i})\cos\kappa\lambda_{T_i} + R^2 C(\theta_{T_i}) \cos 2\kappa\lambda_{T_i}),
		\end{align*}
	\end{linenomath}
	where
	\begin{linenomath}
		\begin{align*}
		A(\theta)&=\tfrac K 2 (2\Omega+K)\cos^2\theta +\tfrac 1 4 K^2 \cos^{2\kappa} \theta \Big( (\kappa+1)\cos^2\theta + (2\kappa^2-\kappa-2)-2\kappa^2\cos^{-2}\theta \Big),
		\\
		B(\theta)&=\frac{2(\Omega+K)K}{(\kappa+1)(\kappa +2)}\cos^{\kappa}\theta \Big( (\kappa^2+2\kappa+2)-(\kappa+1)^2\cos^2 \theta \Big),
		\\
		C(\theta)&= \tfrac{1}{4} K^2 \cos^{2\kappa}\theta \Big( (\kappa+1)\cos^2\theta-(\kappa+2) \Big).
		\end{align*}                                \end{linenomath}		
	For our method, we need the directional magnitude of the velocity which is  $V_{ij}=\mathbf u_{ij} \cdot \mathbf{n}_{ij}$ at edge $ij$.
	\\
	In Fig.~\ref{fig:rossby}, it can be seen that the main features of the evolution of the Rossby--Haurwitz wave solution are reproduced correctly. The time series of the errors in the energy and potential enstrophy are depicted in Fig.~\ref{fig:errorrossby}.

	\section{Long term simulation}\label{sec:LongTermSimulations}
	
	In the previous section, we have reproduced some of the standard test cases proposed in~\cite{wil92} for testing novel numerical schemes for the shallow-water equations. Note that these test cases require the integration of the shallow-water equations for relatively short time intervals, with the longest test case being integrated for $t=15$ days. A main motivation for developing a geometric numerical integrator is that they should be suitable for longer integrations. 
	
	To test the ability of the variational discretization of shallow-water equations to carry out longer integration experiments, we \changed{revisit here} Case 2, the nonlinear geostrophic flow on a rotating sphere and Case 3, the flow over a mountain. We test the spatial variational integrator with two different time discretizations, the (variational) Cayley transform and the (non-variational) standard method that applies a Crank--Nicolson time discretization of the continuity equation. Both test cases are integrated for \changed{a 50-day} period.  The time series of the total energy and potential enstrophy errors are depicted in Fig.~\ref{fig:LongTimeSteadyState} and~\ref{fig:LongTimeMountain}. While both time integrators produce reasonable error time series for the short-term integration \changed{up to, say, 15 days, for long-term simulations} only the variational scheme based on the Cayley transform shows hardly any error trends \changed{for both the nonlinear geostrophic flow and the flow over the mountain test cases.} Notably, while the standard time integrator performs well for Case 2, 
	this non-variational method shows a clear trend to lose energy for Case 3.
	
	Moreover, we note that for all test cases studied, the Cayley transform method shows first order convergence of the energy error with respect to the time step while for the standard method, though showing in general good energy conservation properties, such convergence behaviour cannot be guaranteed. For instance, while the energy error converges at first order for the 
	Rossby--Haurwitz wave case, the standard time integrator yields only zero order convergence for the steady state case (not shown).
	
	
	{\color{black}

		\changed{
			The slight positive drift in energy conservation that is visible in
			the top-left panel of Fig.~\ref{fig:LongTimeSteadyState} is probably
			related to the fact that the presented RSW scheme is not a fully
			variational integrator in time (see section Time discretizations). In
			particular, when using non-uniform grids (such as those used here,
			with both pentagonal and hexagonal cells in the dual mesh), this effect is enhanced. This drift gets smaller with increased resolution, and we observe very good long-term behaviour over a time period of a couple of months, a reasonable time period in atmospheric and ocean modelling. When considering longer integration times of up to one year, for instance, 
			the drift is more enhanced but stays, nevertheless, very small,
			as illustrated in Figure~\ref{fig:LongTimeRegIrreg} (right) for the
			geostrophic flow test case. Similarly to results shown in
			\cite{baue17a}, these trends can be related to the irregularity of the
			mesh.  In particular,  we don't observe any trend in energy error when
			simulating a steady vortex solution on an f-plane approximation of the
			sphere with a uniform mesh (cf. Figure~\ref{fig:LongTimeRegIrreg},
			left), but do see such a trend when using a non-uniform mesh (Figure~\ref{fig:LongTimeRegIrreg}, middle), until the vortex decays. 
			The study of the influence of the grid structures and the construction of a fully variational time integrator to avoid these trends will be part of future work.
		}

		\begin{figure}[!ht]
			\includegraphics[width=\textwidth]{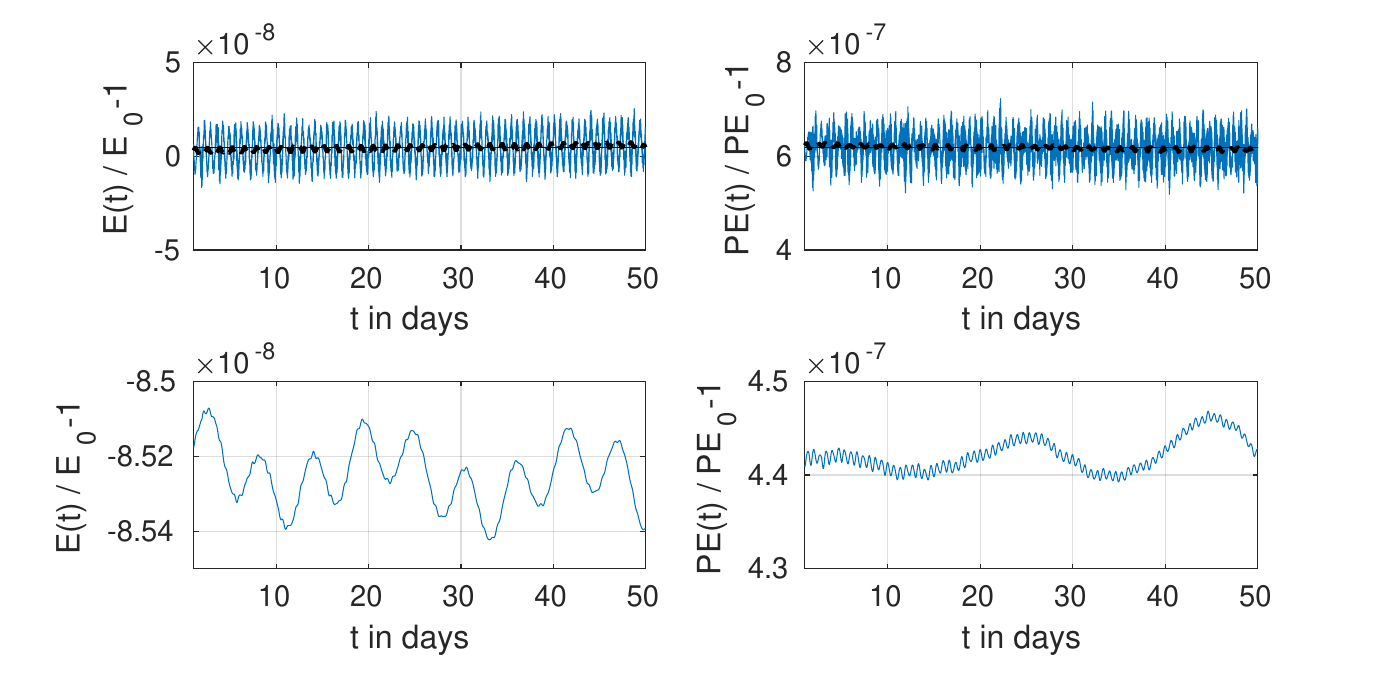}				
			\caption{\footnotesize{Comparison between the Cayley transform (top panels) and standard time integrator (bottom panels) for the nonlinear geostrophic flow (Case 2) over 50 days. The dotted line indicates the moving mean and the solid straight line the overall mean.}}
			\label{fig:LongTimeSteadyState}
		\end{figure}
		
		
		\begin{figure}[!ht]
			\includegraphics[width=\textwidth]{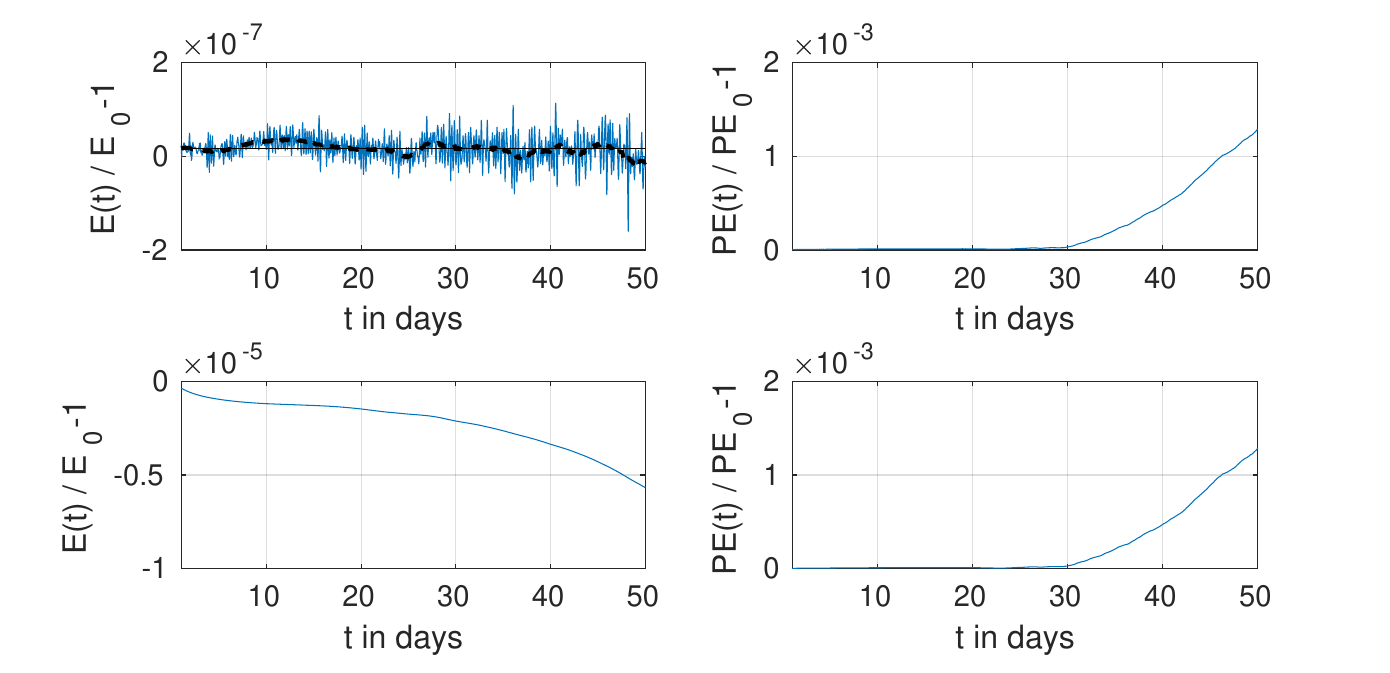}				
			\caption{\footnotesize{Comparison between the Cayley transform (top panels) and standard time integrator (bottom panels) for the flow over the mountain (Case 3) over 50 days.  The dotted line indicates the moving mean and the solid straight line the overall mean.}}
			\label{fig:LongTimeMountain}
		\end{figure}
		
		\begin{figure}[!ht]
			\includegraphics[width=\textwidth]{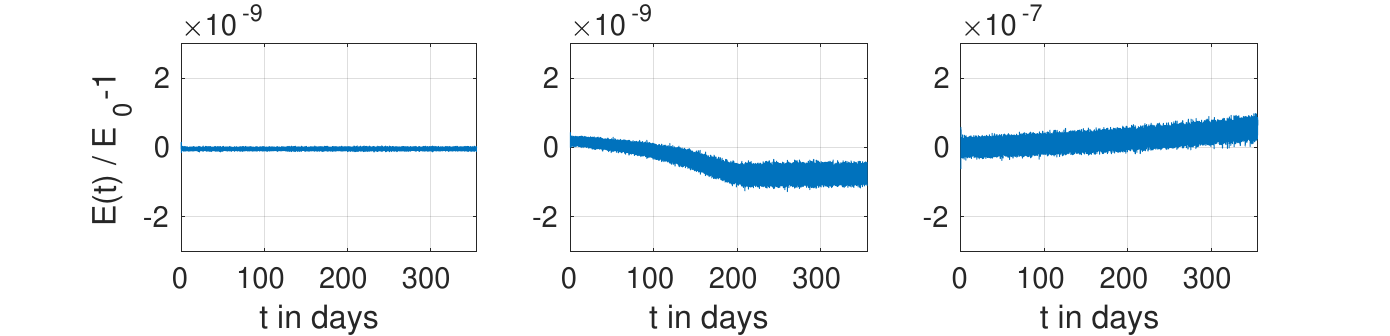}				
			\caption{\changed{\footnotesize{Left and middle plots show the conservation of energy for an isolated vortex (for details see \cite{baue17a}) on a uniform and non-uniform mesh ($64^2$ triangles) on the f-plane over a period of 1 year. The right plot shows the conservation of energy for Case 2 on the sphere (10242 Voronoi cells) for 1 year. }}}
			\label{fig:LongTimeRegIrreg}
		\end{figure}
		
		\section{Conclusions and outlook}\label{sec:ConclusionsVariationalIntegrator}
		
		In this paper, we have constructed a space--time variational discretization for the shallow-water equations on a rotating sphere. This discretization is an extension of the variational integrator proposed in~\cite{baue17b} for the rotating shallow-water equations on the plane. We have carried out some of the standard benchmark tests proposed in~\cite{wil92} and \changed{illustrated that the 
			discretization converges at the order of about 0.5 to 1 (constrained by the convergence order of the 
			divergence operator) to}				
		reference solutions of the shallow-water equations on the rotating sphere. \changed{As expected from the variational integrator, the magnitude of small error fluctuations of the energy around a conserved long-term mean reduce with first order with respect to the time step size.} 
		All numerical tests carried out demonstrate the excellent conservation properties of the variational integrator, regarding the conservation of energy\changed{, mass, and potential vorticity.}
		\changed{Potential enstrophy is not an invariant of the discrete equations by construction; nevertheless 
			it is well preserved.}

		We would like to stress that the variational integrator on the sphere proposed in this paper exactly conserves the lake-at-rest steady state solution over arbitrary bottom topography. In addition, mass is conserved up to machine precision. Both factors are of significant importance in tsunami propagation models, since they prevent the introduction of spurious waves by the discretization scheme. As such, the variational integrator can serve as dynamical core for a general purpose tsunami model, \changed{in particular as our framework permits application of different boundary conditions (e.g., free-slip), cf. \cite{baue17b}.} A tsunami model requires a suitable inundation model to handle the time-dependent wet--dry interface. The development of such an inundation model for the presented variational integrator is currently underway and will be presented elsewhere, as will be the addition of a framework for dynamic grid adaptation, which is another essential ingredient for a tsunami propagation and inundation model.
		
		\section*{Acknowledgments}
		
		This research was undertaken, in part, thanks to funding from the Canada Research Chairs program, the NSERC Discovery Grant program and the InnovateNL CRC LeverageR{\&}D program. WB has received funding from the European Union's Horizon 2020 research and innovative programme under the Marie Sklodowska-Curie grant agreement No 657016. AB is a recipient of an APART Fellowship of the Austrian Academy of Sciences. FGB was partially supported by the ANR project GEOMFLUID, ANR-14-CE23-0002-01. The authors thank Chris Eldred for providing us with the grid generator on the sphere.
		
		\appendix
		\section{Appendix}
		
		\subsection{Differential forms, flat operator, and Lie derivatives}\label{Appendix_forms}
		
		Let us consider a manifold $M$ of dimension $n$, for instance the sphere $\mathcal{S}\subset \mathbb{R}^3$ with $n=2$. For $x\in M$, the vector space $T_xM$ consists of all the (tangent) vectors at the point $x$. A vector field $\mathbf{u}$ on $M$ is a smooth map that associates to each point $x\in M$ a tangent vector $\mathbf{u}(x)\in T_xM$. In local coordinates $(x^1,...,x^n)$ of $M$, a vector field is written $\mathbf{u}(x)=\sum_{i=1}^n \mathbf{u}^i(x)\frac{\partial }{\partial x^i}$.
		
		\medskip
		
		A $k$-form $\omega$ is a skew-symmetric tensor field of rank $k$ on $M$. That is, it is a smooth map that associates to each point $x\in M$ a multi-linear map $\omega(x)$ that takes $k$ tangent vectors at $x$ as input and returns a real number:
		\[
		\omega(x): T_xM\times .... \times T_xM \rightarrow\mathbb{R},\quad (v^1_x,...,v^k_x)\mapsto \omega(x)\big(v^1_x,...,v^k_x \big).
		\]
		A zero-form on $M$ is thus a function on $M$.
		In local coordinates $(x_1,...,x_n)$ of $M$, a $k$-form reads $\omega= \sum_{i_1<...<i_k}\omega_{i_1...i_k}dx^{i_1}\wedge ... \wedge dx^{i_k}$. For instance one-forms and two-forms on $M=\mathbb{R}^2$ read $\omega=\omega_1 dx^1+ \omega_2dx^2$ and $\omega= \omega_{12}dx^1\wedge dx^2$ in coordinates $(x^1,x^2)$.
		The exterior differential is a mapping $\mathbf{d}$ from $k$-forms to $(k+1)$-forms that extends the notion of the differential of a function to differential forms. For instance, in coordinates, the differential of a $k$-form $\omega=\sum_{i_1<...<i_n}\omega_{i_1...i_k}dx^{i_1}\wedge ... \wedge dx^{i_1}$ is the $(k+1)$-form $\mathbf{d}\omega =\sum_{j=1}^n \sum_{i_1<...<i_n}\frac{\partial \omega_{i_1...i_k}}{\partial x^j} dx^j\wedge dx^{i_1}\wedge ... \wedge dx^{i_1}$.
		Given a vector field $\mathbf{u}$ and a $k$-form $\omega$, the inner product of $\mathbf{u}$ with $\omega$ is the $(k-1)$-form $\mathbf{i}_\mathbf{u}\omega$ defined by $\mathbf{i}_\mathbf{u}\omega (x)\big(v_x^1,...,v^{k-1}_x\big)= \omega(x) \big(\mathbf{u}(x),v_x^1,...,v^{k-1}_x\big)$.
		\medskip
		
		A Riemannian metric on $M$ is a symmetric and positive-definite tensor field $\gamma$ of rank $2$ on $M$.  That is, it is a smooth map that associates to each point $x\in M$ a symmetric bilinear map $\gamma(x)$ that takes $2$ tangent vectors as input and returns a real number:
		\[
		\gamma: T_xM\times  T_xM \rightarrow\mathbb{R},\quad (v_x^1,v_x^2)\mapsto \gamma(x)\big(v_x^1,v_x^2\big),
		\]
		such that $\gamma(v_x,v_x)> 0$, for all $v_x\in T_xM$ with $v_x\neq 0$. In coordinates, we write $\gamma=\sum_{i,j=1}^n \gamma_{ij}dx^idx^j$. The flat operator associated to $\gamma$ is the linear map that sends a vector field $\mathbf{u}$ to the 1-form $\mathbf{u}^\flat$ defined by $\mathbf{u}^\flat (x)(v_x)= \gamma(x)(\mathbf{u}(x),v_x)$, for all $v_x\in T_xM$. In local coordinates, if $\mathbf{u}=\sum_{i=1}^n \mathbf{u}^i\frac{\partial }{\partial x^i}$, then $\mathbf{u}^\flat =\sum_{i,j=1}^n \gamma_{ij}\mathbf{u}^j dx^i$. The volume form associated to the Riemannian metric $\gamma$ is the $n$-form $d\sigma$, given locally by $d\sigma = \sqrt{\operatorname{det}\gamma_{ij}}dx^1\wedge ... dx^n$.
		
		\medskip
		
		Given a diffeomorphism $\varphi$ of $M$, the pull-back of a $k$-form $\omega$ is the $k$-form $\varphi^*\omega$ defined by
		\[
		\varphi^*\omega(x)\big(v^1_x,...,v^k_x\big):= \omega(\varphi(x))\big(T_x\varphi(v_x^1),...,T_x\varphi(v_x^k)\big),
		\]
		where $T_x\varphi:T_xM\rightarrow T_{\varphi(x)}M$ is the tangent map (the derivative) of $\varphi$. The Lie derivative of a $k$-form $\omega$ with respect to a vector field $\mathbf{u}$ is the $k$-form $\mathcal{L} _\mathbf{u} \omega$ defined by $\mathcal{L} _\mathbf{u} \omega= \left.\frac{d}{dt}\right|_{t=0}\varphi_t^* \omega$, where $\varphi_t$ is the flow of the vector field $\mathbf{u}$, i.e. $\frac{d}{dt}\varphi_t(x)= \mathbf{u}( \varphi_t(x))$, with $\varphi_{t=0}(x)=x$.
		The divergence of a vector field relative to the Riemannian metric $\gamma$ is the function $\operatorname{div}\mathbf{u}$ defined by $\mathcal{L}_\mathbf{u}d\sigma= (\operatorname{div}\mathbf{u})d\sigma$.
		
		\medskip
		
		For the treatment of fluid mechanics, we also need the notion of Lie derivative of a $k$-form density, denoted $\mathsf{L}_\mathbf{u}$ to distinguish it from the Lie derivative $\mathcal{L}_\mathbf{u}$ of a $k$-form. Let us assume that a Riemannian metric $\gamma$ is fixed (for the sphere $\mathcal{S}$ we chose the Riemannian metric induced from the Euclidean metric on $\mathbb{R}^3$). In this case, $k$-form densities can be identified with $k$-forms. The action of a diffeomorphism $\varphi$ on a $k$-form density $\omega$ is $\omega \bullet\varphi:= (\varphi^*\omega )J\varphi$, where $J\varphi$ is the Jacobian of the diffeomorphism with respect to the Riemannian metric $\gamma$, defined by $\varphi ^*d\sigma= J\varphi d\sigma$. The Lie derivative is defined, similarly as before, by $\mathsf{L} _\mathbf{u} \omega= \left.\frac{d}{dt}\right|_{t=0} \omega\bullet \varphi_t=\left.\frac{d}{dt}\right|_{t=0}(\varphi_t^* \omega)J \varphi_t$, where $\varphi_t$ is the flow of the vector field $\mathbf{u}$. In \S\ref{sec:VariationalPrinciple}, this Lie derivative was applied to the fluid momentum density $\alpha=\frac{\delta \ell}{\delta\mathbf{u}}$ (a one-form density) and to the fluid depth $h$ (a zero-form density), and reads
		\[
		\mathsf{L}_\mathbf{u}\alpha = \mathbf{i}_\mathbf{u} \mathbf{d}\alpha+ \mathbf{d}\mathbf{i}_\mathbf{u}\alpha+\alpha\operatorname{div}\mathbf{u}\qquad\text{and}\qquad \mathsf{L}_\mathbf{u}h= \operatorname{div}(h\mathbf{u}).
		\]
		In local coordinates, we have $(\mathsf{L}_\mathbf{u}\alpha)_i=\sum_{k=1}^n\partial_k (\alpha_i\mathbf{u}^k)+\sum_{k=1}^n\alpha_k\partial_i \mathbf{u}^k$ and $\mathsf{L}_\mathbf{u}h= \sum_{k=1}^n\partial_k (h \mathbf{u}^k)$.
		
		\subsection{Lie group, Lie algebra, and actions}

		Recall that a Lie group $G$ is a manifold and a group, such that the group operations are smooth maps with respect to the manifold structure. Basic examples are the general linear group $\mathsf{GL}(n)=\{A\in \operatorname{Mat}(n)\mid \operatorname{det}A\neq 0\}$ or the special orthogonal group $\mathsf{SO}(n)=\{R \in \operatorname{GL}(n) \mid R^\mathsf{T}R=I_n,\;\;\operatorname{det}R>0\}$, where $I_n$ is the $n\times n$ identity matrix. The tangent space to a Lie group at the identity, denoted $\mathfrak{g}=T_eG$, is called the Lie algebra of the Lie group and is naturally endowed with a Lie bracket $[\,,]$. For matrix Lie groups, the Lie bracket is the usual commutator of matrices $[A,B]=AB-BA$.
		
		\medskip

		Important Lie groups for this paper are the group $\operatorname{Diff}(M)$ of all smooth diffeomorphisms of the manifold $M$ and its discrete version $\mathsf{D}(\mathbb{M})$ defined in \eqref{DD_all}. The Lie algebra of $\operatorname{Diff}(M)$ is the space of vector fields $\mathbf{u}$ on $M$, endowed with (minus) the Lie bracket of vector fields $[\mathbf{u},\mathbf{v}]= \mathbf{u}\cdot \nabla\mathbf{v}-\mathbf{v}\cdot \nabla\mathbf{u}$. The Lie algebra of the Lie group $\mathsf{D}(\mathbb{M})=\left\{q\in \operatorname{GL}(N)^+\mid q\cdot\boldsymbol{1}=\boldsymbol{1}\right\}$ of discrete diffeomorphisms is obtained by taking the $\varepsilon$-derivative of the relation $q_\epsilon \cdot \boldsymbol{1}=\boldsymbol{1}$, for a curve $q_\varepsilon \in \mathsf{D}(\mathbb{M})$ with $q_{\varepsilon=0} =I_N$. We obtain the Lie algebra $\mathfrak{d}(\mathbb{M})=\{A\in \operatorname{Mat}(N)\mid A\cdot \boldsymbol{1}=0\}$ as given in \eqref{Lie_algebra}.
		
		\medskip
		A (right) action of a Lie group $G$ on a manifold $M$ is a map $(g,x)\in G\times M \rightarrow x\cdot g \in M$ such that $x\cdot e=x$ and $x\cdot (gh)= (x\cdot g) \cdot h$, for all $x\in M$ and $g,h\in G$, where $e$ is the neutral element in $G$. The relevant action for the present paper is the action of $\operatorname{Diff}(M)$ on $k$-forms densities $(\varphi, \omega)\mapsto \omega\bullet \varphi=(\varphi^* \omega)J\varphi$. In particular, we have the action on one-form densities $\alpha\bullet \varphi= (\varphi^* \alpha)J\varphi$ and on zero-form densities $h\bullet \varphi= (h\circ \varphi)J\varphi$, where $\alpha$ represents the momentum fluid density and $h$ represents the fluid depth.

		\subsection{Euler--Poincar\'e variational principle}\label{Appendix_EP_principle}

		We quickly review the Euler--Poincar\'e theory by applying it to simpler examples, before considering the rotating shallow water equations.
		Euler--Poincar\'e variational principle applies to dynamical systems whose configuration manifold is a Lie group $G$ and whose dynamics is given by the Euler--Lagrange equations associated to a Lagrangian which is invariant under a subgroup of $G$. We refer to \cite{holm98Ay} for a complete treatment.
		
		\paragraph{The rigid body.} One of the simplest example of dynamical system on a Lie group is a rigid body moving about its center of mass. The configuration Lie group is $\mathsf{SO}(3)$, with $R\in \mathsf{SO}(3)$ describing the orientation of the body. The dynamics of the body is described by the Euler--Lagrange equations for the Lagrangian $L:T\mathsf{SO}(3)\rightarrow\mathbb{R}$ given by the kinetic energy of the body. Euler--Lagrange equations arise from the Hamilton principle
		\begin{equation}\label{HP_body}
		\delta\int_0^TL(R, \dot R){\rm d}t=0,
		\end{equation}
		with respect to variations $\delta R$ with $\delta R(0)=\delta R(T)=0$. This Lagrangian is invariant under the left action of $SO(3)$, i.e., it satisfies $L(AR, A\dot R )= L(R, \dot R )$, for all $A\in\mathsf{SO}(3)$. This invariance allows one to define the reduced Lagrangian $\ell:\mathfrak{so}(3)\rightarrow\mathbb{R}$, with $ \mathfrak{so}(3)$ the Lie algebra of $\mathsf{SO}(3)$, by $L(R, \dot R)=\ell(R^{-1}R)$. The skew-symmetric $3\times 3$ matrix $\Omega=R^{-1}\dot R$ is the angular velocity of the rigid body in the body frame. It is natural to ask if it is possible to write the Hamilton principle \eqref{HP_body} directly in terms of the body angular velocity $\Omega$ without referring to the initial variables $R$ and $\dot R$. This is indeed the case if one considers the constrained variations of $\Omega=R^{-1}\dot R$ induced by the free variations of $R$. A straightforward computation shows that 
		\begin{equation}\label{EP_variations}
		\delta\Omega= \dot A+ [\Omega, A],
		\end{equation}
		where $A= R^{-1}\delta R$ is an arbitrary curve of skew-symmetric $3\times 3$ matrices vanishing at $t=0,T$. The Hamilton principle \eqref{HP_body} is thus equivalent to the variational principle
		\begin{equation}\label{EP_body}
		\delta\int_0^T\ell(\Omega)dt=0,
		\end{equation}
		with respect to variations $\delta \Omega$ of the form \eqref{EP_variations}. This is the Euler--Poincar\'e principle, as it applies to the rigid body. We do not present explicitly the Lagrangian and the equations for the rigid body as we will not use them. We will however use the variational treatment of the rigid body, especially the  passage from \eqref{HP_body} to \eqref{EP_body}, as a guide for the understanding of the more involved treatment of fluids below.
		
		\paragraph{The incompressible fluid.} 
		For the incompressible fluid on a Riemannian manifold $M$, the configuration Lie group is the group $\operatorname{Diff}_{\rm vol}(M)=\{\varphi\in\operatorname{Diff}(M)\mid J\varphi=1\}$ of diffeomorphisms of $M$ that preserve the volume. The equations of motions are given by the Euler--Lagrange equations for the Lagrangian $L:T\operatorname{Diff}_{\rm vol}(M)\rightarrow\mathbb{R}$ given by the kinetic energy of the fluid, namely
		\[
		L(\varphi, \dot \varphi)=\int_M \frac{1}{2} \gamma(\dot\varphi, \dot\varphi)d\sigma,
		\]
		where we recall that $d\sigma = \sqrt{\operatorname{det}\gamma_{ij}}dx^1\wedge ... dx^n$ is the Riemannian volume form associated to $\gamma$. The Hamilton principle reads
		\begin{equation}\label{HP_fluid}
		\delta\int_0^TL(\varphi, \dot \varphi){\rm d}t=0,
		\end{equation}
		with respect to variations $\delta \varphi$ with $\delta \varphi(0)=\delta \varphi(T)=0$. In the same way the Lagrangian of the rigid body was left-invariant, the Lagrangian of the incompressible fluid is right-invariant: $L( \varphi\circ \psi,  \dot \varphi\circ \psi)=  L(\varphi, \dot \varphi)$ for all $\psi\in \operatorname{Diff}_{\rm vol}(M)$. This invariance allows one to define the reduced Lagrangian $\ell:\mathfrak{X}_{\rm vol}(M)\rightarrow\mathbb{R}$, by $L(\varphi, \dot \varphi)=\ell(\dot \varphi\circ\varphi^{-1})$, where $\mathfrak{X}_{\rm vol}(M)$ is the Lie algebra of $\operatorname{Diff}_{\rm vol}(M)$, given by divergence free vector fields on $M$. The reduced Lagrangian is given by 
		\begin{equation}\label{reduced_l_fluid}
		\ell(\mathbf{u})=\frac{1}{2}\int_M\gamma(\mathbf{u},\mathbf{u})d\sigma.
		\end{equation}
		As before, it is natural to ask if it is possible to write the Hamilton principle \eqref{HP_fluid} directly in terms of the Eulerian velocity $\mathbf{u}$ without referring Lagrangian fluid trajectory and velocity $\varphi$ and $\dot \varphi$. This is indeed the case if one considers the constrained variations of $\mathbf{u}= \dot \varphi\circ\varphi^{-1}$ induced by the free variations of $\varphi$. A straightforward computation shows that 
		\begin{equation}\label{EP_fluid}
		\delta\mathbf{u}= \partial _t\mathbf{v}+ [\mathbf{u}, \mathbf{v}],
		\end{equation}
		where $\mathbf{v}= \delta \varphi\circ\varphi^{-1}$ is an arbitrary curve of divergence free vector fields vanishing at $t=0,T$. The Hamilton principle \eqref{HP_fluid} is thus equivalent to the variational principle
		\[
		\delta\int_0^T\ell(\mathbf{u})dt=0,
		\]
		with respect to constrained variations $\delta \mathbf{u}$ of the form \eqref{EP_fluid}. This is the Euler--Poincar\'e principle, as it applies to the incompressible fluid. A direct application of this principle yields
		\begin{equation}\label{application_EP_fluid}
		\delta\int_0^T\ell(\mathbf{u})dt= \int_0^T\int_M \frac{\delta \ell}{\delta \mathbf{u}}\cdot \delta\mathbf{u} d\sigma dt=\int_0^T\int_M\frac{\delta \ell}{\delta \mathbf{u}}\cdot\big( \partial _t\mathbf{v}+ [\mathbf{u}, \mathbf{v}] \big)d\sigma dt= -\int_0^T\int_M\big(\partial_t\frac{\delta \ell}{\delta \mathbf{u}}+ \mathcal{L}_\mathbf{u}\frac{\delta \ell}{\delta \mathbf{u}}\big)\cdot \mathbf{v}d\sigma dt,
		\end{equation}
		where in the first equality the functional derivative $\frac{\delta\ell}{\delta\mathbf{u}}$ is defined as the one-form such that
		\[
		\left.\frac{d}{d\varepsilon}\right|_{\varepsilon=0}\ell(\mathbf{u}+\varepsilon\mathbf{v})=\int_M \frac{\delta \ell}{\delta \mathbf{u} }\cdot   \mathbf{v}\,d\sigma,
		\]
		for arbitrary divergence free vector field $\mathbf{v}$. Since \eqref{application_EP_fluid} has to be zero for all divergence free vector fields $\mathbf{v}$, by the Hodge decomposition there exists a function $q$ such that
		\[
		\partial_t\frac{\delta \ell}{\delta \mathbf{u}}+ \mathcal{L}_\mathbf{u}\frac{\delta \ell}{\delta \mathbf{u}}=-\mathbf{d}q.
		\]
		For the Lagrangian \eqref{reduced_l_fluid} of the incompressible fluid, we have $\frac{\delta \ell}{\delta \mathbf{u}}=\mathbf{u}^\flat$, hence using $\mathcal{L}_\mathbf{u}\mathbf{u}^\flat =\mathbf{i}_\mathbf{u}\mathbf{d}\mathbf{u}^\flat + \mathbf{d}(\mathbf{i}_\mathbf{u}\mathbf{u}^\flat)$, we get the Euler equations in the form $\partial_t\mathbf{u}^\flat + \mathbf{i}_\mathbf{u}\mathbf{d}\mathbf{u}^\flat= - \mathbf{d}p$, with $p= q+ \mathbf{i}_\mathbf{u}\mathbf{u}^\flat$.
		
		\paragraph{The rotating shallow water equations.} The configuration Lie group for the rotating shallow water equations on the sphere is the group $\operatorname{Diff}(\mathcal{S})$ of all diffeomorphisms of $\mathcal{S}$. While the equations are still given by the Euler--Lagrange equations for $L$, obtained from the Hamilton principle $\delta\int_0^TL(\varphi, \dot\varphi)dt=0$, a major difference with the previous two cases is the symmetry group of $L$. It is no longer given by the whole configuration Lie group but by a subgroup, namely, the subgroup $\operatorname{Diff}(\mathcal{S})_{h_0}\subsetneq \operatorname{Diff}(\mathcal{S})$ of diffeomorphisms that preserves the initial fluid depth $h_0$. As a consequence, the reduced Lagrangian not only depends on the Eulerian velocity $\mathbf{u}$ as earlier, but also on the current fluid depth $h$. More precisely, we have
		\begin{equation}\label{def_reduced_RSW}
		L(\varphi, \dot \varphi)=\ell(\mathbf{u},h),\quad \text{for}\quad \mathbf{u}=\dot\varphi\circ\varphi^{-1}\quad\text{and}\quad h=(h_0\circ\varphi^{-1})J\varphi^{-1}.
		\end{equation}
		From these relations, the free variations $\delta\varphi$ of the fluid trajectory $\varphi$ induce the following constrained variations of $\mathbf{u}$ and $h$:
		\[
		\delta\mathbf{u}= \partial _t\mathbf{v}+ [\mathbf{u}, \mathbf{v}]\quad\text{and}\quad\delta h = - \operatorname{div}(h\mathbf{v}),
		\]
		where $\mathbf{v}= \delta \varphi\circ\varphi^{-1}$ is an arbitrary curve of divergence free vector fields vanishing at $t=0,T$. This is the variational principle given in \eqref{EP_principle}.
		We now show the derivation of \eqref{EP_general}. We have
		\begin{equation}\label{application_EP_RSW}
		\begin{aligned}
		\delta\int_0^T\ell(\mathbf{u},h)dt&= \int_0^T \int_M\Big(\frac{\delta \ell}{\delta \mathbf{u}}\cdot \delta\mathbf{u}+\frac{\delta \ell}{\delta h} \delta h\Big) d\sigma dt=\int_0^T\int_M\Big(\frac{\delta \ell}{\delta \mathbf{u}}\cdot\big( \partial _t\mathbf{v}+ [\mathbf{u}, \mathbf{v}]\big) - \frac{\delta \ell}{\delta h} \operatorname{div}(h\mathbf{v})\Big)d\sigma dt\\
		&= \int_0^T\int_M\big(-\partial_t\frac{\delta \ell}{\delta \mathbf{u}}- \mathsf{L}_\mathbf{u}\frac{\delta \ell}{\delta \mathbf{u}}+ h\mathbf{d}\frac{\delta \ell}{\delta h}\big)\cdot \mathbf{v}d\sigma dt,
		\end{aligned}
		\end{equation}
		where we used integration by parts with $\mathbf{v}(0)=\mathbf{v}(T)=0$. Since $\mathbf{v}$ is arbitrary, \eqref{EP_general} follows.
		
		Replacing the functional derivatives \eqref{FD_RSW} for the RSW Lagrangian in \eqref{EP_general}, we get
		\begin{equation}\label{intermediate_EP}
		\partial_t(h(\mathbf{u}^\flat + \mathbf{R}^\flat))+ \mathsf{L}_\mathbf{u} \big(h(\mathbf{u}^\flat + \mathbf{R}^\flat)\big)= h\mathbf{d}\Big(\frac{1}{2}\gamma(\mathbf{u},\mathbf{u}) +\gamma(\mathbf{R},\mathbf{u})- g(h+B) \Big).
		\end{equation}
		Using the formulas for the Lie derivative recalled earlier, we have $\mathsf{L}_\mathbf{u} \big(h(\mathbf{u}^\flat + \mathbf{R}^\flat)\big)=\operatorname{div}(h\mathbf{u})(\mathbf{u}^\flat + \mathbf{R}^\flat)+ \mathbf{i}_\mathbf{u}\mathbf{d}(\mathbf{u}^\flat + \mathbf{R}^\flat)+ \mathbf{d}\big(\mathbf{i}_\mathbf{u} (\mathbf{u}^\flat + \mathbf{R}^\flat)\big)$. Inserting this relation into \eqref{intermediate_EP} and simplifying further by using the advection equation $\partial_th +\operatorname{div}(h\mathbf{u})=0$, which follows from the second equation in \eqref{def_reduced_RSW}, we get the RSW equations in the form given in \eqref{one_form_RSW}.
		

		\color{black}

\footnotesize\setlength{\itemsep}{0ex}

	\end{document}

%% file: grid.tex
	
\def\dotMarkRightAngle[size=#1](#2,#3,#4){%
	\draw ($(#3)!#1!(#2)$) -- 
	($($(#3)!#1!(#2)$)!#1!90:(#2)$) --
	($(#3)!#1!(#4)$);
	\path (#3) --node[circle,fill,inner sep=.5pt]{} ($($(#3)!#1!(#2)$)!#1!90:(#2)$);
}

	\begin{tikzpicture}

\coordinate (A) at (0,0);
\coordinate (B) at (6,0);
\coordinate (C) at (3,5.2);
\coordinate (D) at (9,5.2);

\coordinate (Ti) at (3,1.73);
\coordinate (Tj) at (6,3.46);
\coordinate (Tim) at (0,3.46);
\coordinate (Tip) at (3,-1.73);
\coordinate (Tjm) at (6,6.693);
\coordinate (Tjp) at (9,1.73);

\coordinate (Tn1) at (0,6.93);
\coordinate (Tn2) at (3,8.66);

\coordinate [label={[black]left:$f_{ij}$}] (fij) at (4.5,2.6);
\coordinate [label={[blue]left:$h_{ii_-}$}] (hiim) at (1.5, 2.6);
\coordinate [label={[blue]below left:$h_{ii_+}$}] (hiip) at (3, 0);
\coordinate [label={[blue]above right:$h_{jj_-}$}] (hjjm) at (6, 5.2);
\coordinate [label={[blue]right:$h_{jj_+}$}] (hjjp) at (7.5, 2.6);

\draw[gray,pattern=north east lines,pattern color=gray!30] (hiim) -- (C) -- (fij) -- (Ti) -- (hiim);
\draw[gray,pattern=north west lines,pattern color=gray!30] (Ti) -- (fij) -- (B) -- (hiip) -- (Ti);
\draw[gray,pattern=north west lines,pattern color=gray!30] (Tj) -- (fij) -- (C) -- (hjjm) -- (Tj);
\draw[gray,pattern=north east lines,pattern color=gray!30] (Tj) -- (fij) -- (B) -- (hjjp) -- (Tj);

\coordinate [label={[black]below:$K_{i_-}$}] (Kim) at (3,3.46);
\coordinate [label={[black]left:$K_{i_+}$}] (Kip) at (4.5,0.87);
\coordinate [label={[black]below:$K_{j_-}$}] (Kjm) at (4.5,4.33);
\coordinate [label={[black]above:$K_{j_+}$}] (Kjp) at (6,1.73);

\coordinate [label={[black]left:$f_{ij}$}] (fij) at (4.5,2.6);

\dotMarkRightAngle[size=9pt](Ti,fij,B)

\draw[draw, very thick] (A) -- (B) -- (C)--(A);
\draw[draw, very thick] (B) -- (C) -- (D)--(B);

\draw[dashed,blue,very thick] (Tim) -- (Ti);
\draw[dashed,blue,very thick] (Ti) -- (Tip);
\draw[dashed,blue,very thick] (Ti) -- (Tj);
\draw[dashed,blue,very thick] (Tj) -- (Tjm);
\draw[dashed,blue,very thick] (Tj) -- (Tjp);

\draw[dashed,blue] (Tim) -- (Tn1);
\draw[dashed,blue] (Tn1) -- (Tn2);
\draw[dashed,blue] (Tn2) -- (Tjm);

\filldraw[blue] (B) circle  (2pt) node[anchor=north] {$\zeta_+$};
\filldraw[blue] (C) circle  (2pt) node[anchor=south] {$\zeta_-$};
\filldraw[black] (Ti) circle  (2pt) node[anchor=south] {$T_i$};
\filldraw[black] (Tj) circle  (2pt) node[anchor=north] {$T_j$};
\filldraw[black] (Tim) circle (2pt) node[anchor=south] {$T_{i_-}$};
\filldraw[black] (Tip) circle (2pt) node[anchor=east] {$T_{i_+}$};
\filldraw[black] (Tjp) circle (2pt) node[anchor=south] {$T_{j_+}$};
\filldraw[black] (Tjm) circle (2pt) node[anchor=south] {$T_{j_-}$};

	\end{tikzpicture}